\documentclass[11pt,letterpaper,reqno]{amsart}
\usepackage[T1]{fontenc}
\usepackage[latin9]{inputenc}
\usepackage[margin=1in,footskip=.25in]{geometry}
\usepackage{subcaption}
\usepackage{wrapfig}

\usepackage{multicol}
\usepackage{mathtools}
\usepackage{enumerate}
\usepackage{graphicx}
\usepackage{soul}
\usepackage{xcolor}
\usepackage{amssymb}
\usepackage{bigints}
\usepackage{amsmath}
\usepackage{breqn}

\def\c{\mathfrak{c}}

\def\R{\mathbb{R}}

\def\c{c_k}

\usepackage[T1]{fontenc}
\usepackage{makecell}

\theoremstyle{definition}

\numberwithin{equation}{section}
\allowdisplaybreaks
\usepackage{verbatim}
\usepackage{placeins}
\usepackage{cite}
\usepackage{caption}
\usepackage{stmaryrd}
\usepackage{enumerate}
\usepackage{afterpage}
\usepackage{enumitem}
\usepackage{bmpsize}
\usepackage{hyperref}
\usepackage{tabu}
\usepackage{enumitem}
\usepackage{tikz}
\usetikzlibrary{matrix,graphs,arrows,positioning,calc,decorations.markings,shapes.symbols}
\definecolor{dullmagenta}{rgb}{0.4,0,0.4}   
\definecolor{darkblue}{rgb}{0,0,0.4}

\mathtoolsset{showonlyrefs}

\begin{document}
\pagestyle{plain}
\title[On the stationary measure for the open KPZ equation]{Some recent progress on the stationary measure for the open KPZ equation}
\author{Ivan Corwin}

\address[Ivan Corwin]{\parbox{\linewidth}{ Department of Mathematics, Columbia University,
 New York, NY, USA \\ E-mail: ic2354@columbia.edu}}

\begin{abstract}
This note, dedicated in Harold Widom's memory, is an expanded version of a lecture I gave in fall 2021 at the MSRI program ``Universality and Integrability in Random Matrices and Interacting Particle Systems''. I will focus on the behavior of the stationary measure for the open KPZ equation, a paradigmatic model for interface growth in contact with boundaries. Much of this will review elements of my joint work with A. Knizel as well as with H. Shen, as well as subsequent works of W. Bryc, A. Kuznetsov, Y. Wang, and J. Weso{\l}owski and of G. Barraquand and P. Le Doussal. The basis for this advance is fundamental work of B. Derrida, M. Evans, V. Hakim and V. Pasquier from 1993, of T. Sasamoto, M. Uchiyama and M. Wadati from 2003, and of W. Bryc and J. Weso{\l}owski from 2010 and 2017. I will try to explain how all of this fits together, without laboring details for the sake of exposition.

Though this work does not directly follow from Harold Widom's own work, it (and a great deal of my research) is very much inspired by his and Craig Tracy's work on ASEP.
\end{abstract}
\maketitle

\section{Preface}\label{sec:intro}
The purpose of this note is to describe the stationary measure for the open KPZ (Kardar-Parisi-Zhang) equation and some of the ideas related to its construction. The KPZ equation has been a subject of intense study due to its value as a model for random growth and as a singular stochastic PDE, its remarkable connection with integrable systems, its ubiquitous occurrence in seemingly unrelated problems in mathematics and physics, e.g., see the reviews \cite{doi:10.1142/S2010326311300014,HHTak,Quastelspohn,CorwinAMS,Ganguly}

The open KPZ equation is supposed to model stochastic growth in contact with boundaries, or equivalently (through differentiation) stochastic transport between reservoirs. As is often the case in the study of Markov processes, a fundamental question here is to understand the structure of its stationary state -- for instance, whether it is unique and how it can be characterized.

Remarkably, as I will explain below, this stationary measure can be described in terms of the Brownian motion measure reweighted in terms of certain exponential functionals -- a subject which has attracted great interest within other realms of probability theory in its own right, e.g. \cite{comtet_monthus_yor_1998,HariyaYor,10.1214/154957805100000159}.

\subsubsection*{Acknowledgements}
This note is based on a talk I delivered at MSRI in the fall 2021 introductory workshop. My participation in that program was partially supported by the National Science Foundation under Grant No. 1440140, the Chern Professorship from MSRI, and a Visiting Miller Professorship from the Miller Institute for Basic Research in Science. I also wish to acknowledge ongoing support from the NSF through DMS:1811143 and DMS:1937254, the Simon Foundation through a Simons Fellowship in Mathematics (Grant No. 817655), the W.M. Keck Foundation Science and Engineering Grant program, and the Packard Foundation Fellowship for Science and Engineering. I also wish to thank Alisa Knizel and Hao Shen for their collaboration with me related to the content of this note, and to Guillaume Barraquand, Wlodek Bryc, Tomohiro Sasamoto, Yizao Wang, and Lauren Williams for discussions in the past relevant to this note. I first learned about Askey-Wilson processes in 2010 at the conference ``Orthogonal polynomials, applications in statistics and stochastic processes'' held at Warwick, and I wish to thank the organizers for that very valuable conference.
I also appreciate Estelle Basor's efforts in putting together this volume in memory of Harold.
\section{An aside on $q$-Pochhammer symbol asymptotics}\label{Sec:qaside}
I start with something that may seem far from the advertised subject. A devoted reader who gets to the very end of this note (see Section \ref{sec:fullcircle}) will see that, indeed, this is quite key in the asymptotic analysis that leads to the main result presented in this note. Indeed, it is very much in the spirit of Harold's work that in the end, things boil down to involved asymptotic analysis.

The $q$-Pochhammer symbol is defined for $a\in \mathbb{C}$ and $|q|<1$ by
$$
(a;q)_{\infty} := (1-a)(1-qa)(1-q^2 a)\cdots.
$$
This convergent infinite product defines an analytic function (in $a$) which arises in many contexts.
In combinatorics it encodes certain generating functions. Two of the simplest examples are
$$
(-q;q)_{\infty} = \sum_{\lambda \textrm{ strict}}q^{|\lambda|}\qquad \textrm{and} \qquad(q;q)_{\infty}^{-1} = \sum_{\lambda}q^{|\lambda|}
$$
where $\lambda$ denotes a partition (i.e., a weakly decreasing sequence of non-negative integers $\lambda_1\geq \lambda_2\geq \cdots$), strict means  that $\lambda_1>\lambda_2>\cdots$, and $|\lambda|:=\lambda_1+\lambda_2+\cdots$.

The $q$-Pochhammer symbol is key to defining various $q$-deformed variants of classical special functions. These come up, for example, in the Askey-Wilson scheme of orthogonal polynomials, see, for example, Section \ref{sec:AWPP} below or \cite{AskeyWilson}.

The $q$-gamma function is a deformed special function, defined for $z\in \mathbb{C}$ and $|q|<1$ as
$$
\Gamma_q(z) := (1-q)^{1-z} \frac{(q;q)_{\infty}}{(q^z;q)_{\infty}}.
$$
It is not too hard to see that
$$
\lim_{q\to 1} \Gamma_q(z) = \Gamma(z),
$$
for instance by observing convergence to the Euler product formula for the gamma function.

\subsection{Asymptotics}
It is an interesting and (as will be clear much later) valuable question to determine the nature of this convergence. For instance, the gamma function has various asymptotic behaviors as $|z|\to \infty$, in a manner depending on the direction to which $z$ goes to complex infinity. Does the $q$-gamma function enjoy similar asymptotic behavior and is this uniform, in any sense, as $q$ tends to $1$?

In 1984, Moak \cite{Moak} first took up this question, proving a Stirling's type expansion for $\Gamma_q(z)$ with $q\in (0,1)$ fixed and $z$ tending to infinity with $\arg(z)\in(-\pi/2 +\delta, \pi/2 - \delta)$ for any fixed $\epsilon$. Soon after, McIntosh \cite{McIntosh} addressed the case where $z$ is fixed but $q$ tends to $1$ from below, see also work of Daalhuis \cite{Daalhuis} and more recently Katsurada \cite{852f802e715846339a94147fc8f7f3bb} and Zhang \cite{Z}.

The most complicated element in $\Gamma_q(z)$ is the denominator $(q^z;q)_{\infty}$. Consider $z\in\mathbb{C}$ and
$$q=e^{-\epsilon}$$
as $\epsilon$ tends to $0$. In that case, $q^z$ approaches $1$ (for $z$ fixed) and hence many terms in the product defining $(q^z;q)_{\infty}$ will approach 0. In fact, a quick calculation shows that there are order $\epsilon^{-1}$ such terms. This suggests  that $\log (q^z;q)_{\infty}$ will decay like a negative constant times $\epsilon^{-1}$. This is the case, and the constant is remarkably given by $\zeta(2)$. The expansion continues like
$$
\log(q^z;q)_{\infty} = -\tfrac{\pi^2}{6} \epsilon^{-1} -\left(z-\tfrac{1}{2}\right)\log(\epsilon) -\log\left(\tfrac{\Gamma(z)}{2\pi}\right) + \cdots.
$$
The $\Gamma(z)$ term in this expansion is what ultimately leads to the limit $\Gamma_q(z)\to \Gamma(z)$. The $\cdots$ lower order terms hide a lot here. For instance, it is possible to go out to arbitrary order $m\in \mathbb{Z}_{\geq 1}$ in $\epsilon$ so that the lower order terms take the form
$$
\cdots = -\sum_{n=1}^{m} \frac{B_{n+1}(z)B_n}{n(n+1)!}\epsilon^{n} + \textrm{Error}_{m}(\epsilon,z)
$$
where $B_k(z)$ and $B_k$ are the Bernoulli polynomials and Bernoulli numbers, respectively, and the remaining error in this approximation is denoted as $\textrm{Error}_{m}(\epsilon,z)$.

Of course, everything is now moved into studying the behavior of the residual error term $\textrm{Error}_{m}(\epsilon,z)$. Namely, how does it behave as $\epsilon$ and $z$ vary? For $\epsilon$ fixed and $z$ varying in certain regions of the complex plane, or for $\epsilon$ tending to $1$ and $z$ in fixed compact regions of the complex plane, this was understood in the earlier mentioned works. However, what if $z$ varies in a region that is not bounded as $\epsilon$ tends to $0$? It turns out that exactly this type of control is pivotal in my derivation of the open KPZ equation stationary measure with Knizel in \cite{CK} since things eventually boil down to studying asymptotics of integrals involving factors like $(q^z;q)_{\infty}$ for $z$ varying over regions that grow as $q$ goes to $1$. In order to apply the dominated convergence theorem, some type of uniform control is needed. In particular, in my work with Knizel we make an estimate \cite[Proposition 2.2]{CK} that shows that for any $\delta\in (0,1/2)$ and $b\in (m-1,m)$, there exists $C,\epsilon_0>0$ such that for all $\epsilon\in (0,\epsilon_0)$ and all $z\in \mathbb{C}$ with $|\textrm{Im}(z)|<5/\epsilon$,
$$
\left|\textrm{Error}_{m}(\epsilon,z)\right| \leq C \left(\epsilon(1+|z|)^2 + \epsilon^{b} (1+|z|)^{1+2b+\delta}\right).
$$
While we made no claim as to whether this is an optimal bound, it does suffice for our purposes in applying dominated convergence. There is a similar sort of expansion for $(-q^z;q)_{\infty}$, though I will not record it here.

The proof of this asymptotic expansion and error bound in \cite[Proposition 2.2]{CK} builds on work of Zhang \cite{Z}, correcting some mistakes therein and extending from compact domains for $z$ to unbounded ones. The starting point is the following Mellin-Barnes integral representation
$$
\log(q^z;q)_{\infty} = -\frac{1}{2\pi} \int_{c-\iota \infty}^{c+\iota \infty} \Gamma(s)\zeta(s+1)\zeta(s,z) \epsilon^{-s} ds
$$
where  $q=e^{-\epsilon}$, and where $c>1$, $\textrm{Re}(z)>0$, the integral is over a vertical contour, and $\zeta(s,z)$ is the Hurwitz zeta function. Asymptotic of this integral formula requires a well-controlled understanding of asymptotics of the terms in the integrand. For $\Gamma(s)$ and $\zeta(s+1)$ such control is well-known, while for $\zeta(s,z)$ it needs to be derived, again based on Mellin-Barnes integral formulas which express $\zeta(s,z)$ in terms of simpler functions (namely, the gamma and zeta function). Finally, observe that the above Mellin-Barnes formula for $\log(q^z;q)_{\infty}$ was restricted to $\textrm{Re}(z)>0$. The move into the other half of the complex plane is facilitated  by certain identities involving the Jacobi theta functions. I will not give any further details here. An interested reader can find more in the final section of \cite{CK}.

\section{What is the open KPZ equation?}
Since this is not a survey on the KPZ equation as a stochastic PDE, I will try to stay brief on this side of the story and instead refer an interested reader to my earlier survey \cite{doi:10.1142/S2010326311300014} and more recent expository piece with Shen \cite{CorwinShenBAMS}.

\begin{figure}[t]
	\captionsetup{width=.8\linewidth}
	\begin{center}
	\includegraphics[width=3in]{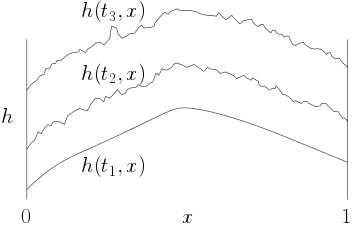}
	\end{center}
	\caption{An illustration of the open KPZ equation evolution at three times $t_1,t_2$ and $t_3$. The height is plotted as a function of space $x$. The height functions, in general, can overlap since they may both grow and shrink. At the left boundary, the slope is maintained as $u$ and at the right it is maintained as $-v$. The initially smooth height profile is roughened with time, as it approaches the stationary measure.}
	\label{fig:KPZ}
\end{figure}

The open KPZ equation is a singular stochastic PDE describing the evolution of a {\it height function} $h(t,x)$ taking real values, with time $t\in \mathbb{R}_{\geq 0}$, and space $x\in [0,1]$. The evolution is a function of a time-space white noise, denoted by $\xi(t,x)$, and described by the equation
$$
\partial_t h = \tfrac{1}{2} \partial_{xx}h + \tfrac{1}{2}(\partial_x h)^2 + \xi
$$
with the boundary conditions for all $t\in\mathbb{R}_{>0}$ depending on two parameters $u,v\in \mathbb{R}$ given by
$$
\partial_x h = \begin{cases} u,&x=0\\-v&x=1.\end{cases}
$$
Since the noise driving this equation is uncorrelated in time, it follows that this defines a Markov process.
At a physical level, this equation is supposed to model the evolution of an interface that is subject to smoothing (the Laplacian), growth normal to the interface (the non-linearity) and a stochastic driving force (the white noise). The boundary condition is suppose to indicate that local to the end of the interval, the slope is maintained at specific values. From a physical perspective, this could be caused by an interaction between the boundary and growing media. See Figure \ref{fig:KPZ} for an illustration of this growth process over time.

Another physical interpretation of this equation is in terms of the dynamics on $\partial_x h$. The resulting stochastic PDE is known as the stochastic Burger's equation and it is supposed to model stochastic conservative transport along the interval $[0,1]$. The boundary conditions then correspond to maintaining specific densities or potentials for reservoirs at the two endpoints of the interval.

When I introduce the open ASEP height and particle process later and describe its limit to the open KPZ equation, both of these interpretations may come into further focus.

\subsection{Long-time and stationary behavior}
The driving question behind studying the open KPZ stationary measure is to understand the long-time behavior of the open KPZ equation. Does the height function have a limit (in distribution) as time goes to infinity? The easy answer to this question is NO. There is nothing pinning the height function down and hence over time the height will drift away to infinity. Of course, the relevant question is not about the absolute height, but rather the relative height -- that is the height function centered by its value at $x=0$. Now there is good reason to believe that the law of this profile should converge as time grows to a limit distribution on height functions which take value 0 at $x=0$, and that the law should not depend on the initial distribution of the profile. I will return to this, as of yet, open problem a bit later.

The law on a height function $h_{u,v}:[0,1]\to \mathbb{R}$ with $h_{u,v}(0)=0$ will be said to be a {\it stationary measure} for the open KPZ equation if when started with $h(0,\cdot) \overset{(d)}{=} h_{u,v}(\cdot)$, it follows that $h(t,\cdot)-h(t,0) \overset{(d)}{=} h_{u,v}(\cdot)$ for all subsequent times $t>0$. Of course, it is reasonable to expect that the long-time limit profile will converge to such a stationary measure, and in the case of the open KPZ equation that for each choice of $u,v$ there is just one such stationary measure. Though these are great questions that should be addressed, here I will focus on the question of how to construct and describe a stationary measure for the open KPZ equation.

\subsection{Defining the open KPZ equation}
\label{sec:defopenKPZ}

Before proceeding in this direction, let me briefly address the highly non-trivial task of  making mathematical sense of what it means to solve the open KPZ equation. Even a n\"aive understanding of time-space white noise suggests that the spatial trajectory of $h(t,x)$ should have the regularity of Brownian motion and hence be H\"older $1/2-$. Thus, making sense of the KPZ equation non-linearity becomes challenging in this case, as does making sense of the meaning of the boundary condition. This can all be done by appealing to Hairer's regularity structure framework, see \cite{GH19}, by smoothing $\xi$ and then performing suitable renormalization (in fact if the smoothing is only on space, Bertini and Cancrini's  earlier approach from \cite{BertiniCancrini} should suffice).

The simplest way, however, to define the open KPZ equation is through the open stochastic heat equation (SHE). For the sake of completeness I will recall this definition from \cite{CS} here. In fact, most of what I discuss in this note can be appreciated without a deep understanding of the stochastic PDE side of this story.

The {\it Hopf-Cole solution} to the open KPZ equation is defined as  $h(t,x) :=\log z(t,x)$ where $z(t,x)$ solves the SHE
$$
\partial_t z = \tfrac{1}{2}\partial_{xx} z + \xi z
$$
with boundary conditions
$$
\partial_x z = \begin{cases} (u-\tfrac{1}{2})z,&x=0,\\-(v-\tfrac{1}{2}) z &x=1.\end{cases}
$$
The inclusion of this $1/2$ factor is just a convention that ensures that the point $u,v=0$ is special (in terms of the phase diagram for the open KPZ equation). The above form of the SHE can be defined and solved via classical methods of semi-linear stochastic PDEs. In particular, a process $z(t,x)$ is a (mild) solution to the SHE if $z(t,\cdot)$ is adapted to the filtration generated by the initial data $z_0(\cdot)$ and the time-space white noise $\xi$ up to time $t$, and if $z$ satisfies for all $t\in \mathbb{R}_{>0}$ and $x\in [0,1]$ the equation
$$
z(t,x) = \int_0^1 p_{u,v}(t;x,y) z_0(y)dy + \int_0^1 \int_{0}^{t} p_{u,v}(t-s;x,y)z(s,y)\xi(ds,dy).
$$
The last integral is in the sense of It\^{o} and $p_{u,v}(t;x,y)$ is the heat kernel on $[0,1]$ satisfying
$$
\partial_t p_{u,v}(t;x,y) = \partial_{xx} p_{u,v}(t;x,y),\qquad p_{u,v}(0;x,y) = \delta_{x=y}
$$
with boundary conditions for all $y\in [0,1]$ and $t\in \mathbb{R}_{>0}$
$$
\partial_x p_{u,v} = \begin{cases} (u-\tfrac{1}{2})p_{u,v} ,&x=0,\\-(v-\tfrac{1}{2}) p_{u,v}&x=1.\end{cases}
$$

The open SHE admits a chaos series expansion and can (though it has not been precisely given in the literature) also be interpreted as a partition function for a continuum directed random polymer model in which the underlying path measure is that of Brownian motion which either dies or splits at the boundaries 0 and 1, depending on the signs of $u$ and $v$.

\section{Constructing the stationary measure}\label{sec:const}
I will relate here the main result of my joint work with A. Knizel, given there as \cite[Theorem 1.2]{CK}. Essentially, we show that for each pair $(u,v)\in \mathbb{R}$ of  boundary parameters there exists a stationary measure for the corresponding version of the open KPZ equation; for general $(u,v)$ we provide some properties of our measure while for $(u,v)$ such that $u+v\geq 0$ we are able to completely characterize the measure in terms of an exact formula for its multi-point Laplace transform. In phrasing our results, there is a bit of subtlety stemming from the fact that we did not show that these are the unique stationary measures (though we conjecture this to be the case).

As suggested above, our first result is that for each pair $(u,v)\in \mathbb{R}$ there exists a stationary measure for the corresponding open KPZ equation with those boundary parameters. We denote a random function distributed according to this measure by $h_{u,v}(\cdot)$. As I will explain later in this note, these measures arise as subsequential limits of a stationary measure for a discrete approximation to the open KPZ equation, namely the open ASEP.

In the special case that $u+v=0$, the random function $h_{u,v}(\cdot)$ has the law of a standard Brownian motion of drift $u=-v$. In fact, in that case the ASEP approximation scheme is not just tight but has a unique limit point. That is not to say that this implies that there is a unique stationary measure in this case, only that the ASEP stationary measures has a unique limit, denoted by $h_{u,-u}(\cdot)$.

For general $(u,v)\in \mathbb{R}$, the increments of $h_{u,v}(\cdot)$ satisfy a property that we call {\it stochastic sandwiching}. The simplest and most useful implication of this sandwiching is the following. There exists a coupling (i.e., a common probability space which supports these random processes) of $h_{u,v}(\cdot)$,  $h_{u,-u}(\cdot)$ and $h_{-v,v}(\cdot)$ such that for all $0\leq x\leq y\leq 1$, almost surely
$$
h_{-v,v}(y)-h_{-v,v}(x) \leq h_{u,v}(y)-h_{u,v}(x)\leq h_{u,-u}(y)-h_{u,-u}(x)
$$
when $u+v> 0$; when $u+v<0$, the above holds with the inequalities swapped.

For $u+v>0$, we completely characterize the law of $h_{u,v}(\cdot)$ in terms of its multi-point Laplace transform. Before stating that, let me relate the simplest version of this formula which characterizes $h_{u,v}(1)$. Physically, this represents the total displacement of the open KPZ height interface on the interval. Due to its Brownianity, in the case where $u+v=0$ this increment is clearly a Gaussian random variable of mean $u=-v$ and variance $1$.
When $u,v>0$ the law of $h_{u,v}(1)$ is determined by the following formula
\begin{equation}\label{eq:netchange}
\mathbb{E}\left[e^{-s h_{u,v}(1)}\right] = e^{s^2/4} \frac{\displaystyle\int_0^{\infty}\!\!\! e^{-r^2} \mu_{s}(r)dr}{\displaystyle\int_0^{\infty}\!\!\! e^{-r^2} \mu_{0}(r)dr},\quad \textrm{for}\quad \mu_{s}(r) = \frac{\left|\Gamma(\tfrac{s}{2}+u+\iota r)\Gamma(-\tfrac{s}{2}+v+\iota r)\right|^2}{\left|\Gamma(2\iota r)\right|^2}
\end{equation}
and where the Laplace variable $s$ is allowed to vary in $(0, 2v)$. When $u+v>0$ yet one of the variables is negative, there is a similar albeit slightly more complicated formula involving a measure with and atomic part as well as an absolutely continuous part.

Turning to the general formula, let $C_{u,v}= 2\mathbf{1}_{u\notin (0,1)} + 2u \mathbf{1}_{u\in (0,1)}$ and consider any $d\in \mathbb{Z}_{\geq 1}$, any $0=x_0<x_1<\cdots<x_d\leq x_{d+1}=1$ and any $c_1,\ldots, c_d$ such that $c_1+\cdots +c_d<C_{u,v}$. Then, letting  $s_k=  c_k+\cdots +c_d$ for $k=1,\ldots, d$ and $s_{d+1}=0$, the multi-point Laplace transform is given by the formula
\begin{equation}\label{eq:laplace}
\mathbb{E}\bigg[e^{-\sum\limits_{k=1}^{d} c_k h_{u,v}(x_k)}\bigg] = \frac{\mathbb{E}\bigg[ e^{\frac{1}{4} \sum\limits_{k=1}^{d+1} (s_k^2-T_{s_k})(x_k-x_{k-1})}\bigg]}{\mathbb{E}\left[e^{-\frac{1}{4} T_{0}}\right]}.
\end{equation}
The expectation on the left side above is over the open KPZ stationary measure on $h_{u,v}(\cdot)$ while on the right it is over the measure on a stochastic process $(T_{s})_{s\in [0,C_{u,v})}$ that we termed the {\it continuous dual Hahn process}. The transition probabilities for this process are given in terms of the orthogonality probability measure for the continuous dual Hahn polynomials, and the measure according to which we initialize this process in the above identity takes a similar form in terms of a ratio of products of gamma functions. Since this involves a number of formulas that probably will not be illuminating, I will not write anything further on this definition here. An interested reader is referred to \cite{CK} or the subsequent work of Bryc \cite{BRYC2022185} for more on this process and the formulas used to define it.

The above characterization of the stationary measure $h_{u,v}(\cdot)$ in the case where $u+v>0$ is, in my own opinion, a bit less than fulfilling. While it uniquely characterizes the law it does not provide a direct description of $h_{u,v}(\cdot)$ in terms of a ``nice'' stochastic process, such as is the case when $u+v=0$ (when $h_{u,-u}(\cdot)$ is Brownian). In fact, it is challenging (though certainly possible) to take the limit of $u+v\searrow 0$ in our Laplace transform formula in order to recover the Laplace transform of the Brownian motion $h_{u,-u}(\cdot)$.

Thankfully (for us), this lack of a nice stochastic process description was fairly quickly resolved, as I explain now.

\section{Inverting the multipoint Laplace transform formula}
It took a few years for Knizel and me to complete our work, proving the above results. Pretty early on, though, we had a good sense of what the final formula should look like and we communicated this to Bryc when he visited us as Columbia for a seminal. He had worked with Wang in \cite{BW} on inverting some similar (albeit less complicated) Laplace transform formulas arising for the limit of the ASEP stationary measure, so it seems reasonable to see if he had ideas on how to proceed. The seed was planted and as luck would have it, around the time that Knizel and I were finishing our paper, Bryc and his collaborators (Kuznetsov, Wang, and Weso{\l}owski) had figured out how to do the desired inversion. In fact, this was extremely helpful since it provided a check for our formulas and helped reveal a few missing factors that we hunted down and repaired.

\begin{figure}[t]
	\captionsetup{width=.8\linewidth}
	\begin{center}
	\includegraphics[width=5in]{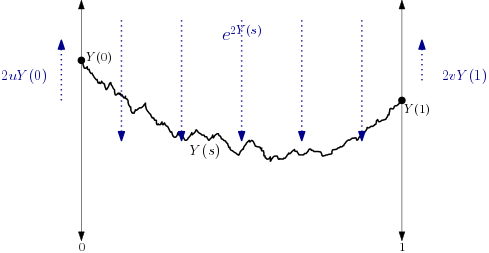}
	\end{center}
	\caption{A depiction of $Y(\cdot)$. In this illustration we assume $u,v>0$ so the starting and ending points are probabilistically rewarded for going up with a linear energy contribution. However, the exponential energy term $e^{Y(s)}$ pushes the curves down. It is this balance of energetic constraints (depicted in dark blue) that results in the reweighted measure being normalizable despite the fact that the free starting point Brownian motion is an infinite measure.}
	\label{fig:huv}
\end{figure}

The description that they arrived at in \cite{BKWW} can be written as follows. The process $x\mapsto h_{u,v}(x)$ for $x\in [0,1]$ has the same law as the process
\begin{equation}\label{BBKWdescription}
x\mapsto  B(x) + Y(0) - Y(x).
\end{equation}
Here $B$  denotes a Brownian motion started at zero and with variance $1/2$ and mean $0$ at time $1$. The process $Y$ is independent of $B$ and has a Radon-Nikodym derivative against a free starting point Brownian motion with variance $1/2$ and no drift which is proportional to (see also Figure \ref{fig:huv})
\begin{equation}\label{BBKWdescriptionRND}
\exp\left(2 u Y(0) - \int\limits_0^{1} e^{2 Y(s)}ds + 2v Y(1)\right).
\end{equation}
By a free starting point Brownian motion, I mean the infinite measure on paths where the starting point is distributed according to Lebesgue measure on $\mathbb{R}$ and the trajectory from there is distributed according to a Brownian motion with that starting point. The time interval for this Brownian motion is $[0,1]$ and the variance at time $1$ is $1/2$. It is a non-trivial calculation to show that despite the free Brownian motion being an infinite measure, after multiplying by the above Radon-Nikodyn derivative, the resulting measure has finite mass and can be normalized to define a probability measure, i.e., the measure on $Y$.  The work in \cite[Proposition 1.6]{BKWW} is valid for $u+v>0$ and $\min(u,v)>-1$, though the description seems to make sense for general $u+v>0$ and that extension should be possible with what is already known.

Soon after the work of \cite{BKWW} was posted, Barraquand and Le Doussal \cite{BLD} used a different set of tools from physics (coming from the study of 1d Liouville quantum gravity and not addressed there with mathematical rigor) to independently invert the Laplace transform. They arrived at a slightly different description that, in fact, has some advantages. They show that the process $x\mapsto h_{u,v}(x)$ for $x\in [0,1]$ has the same law as the process
\begin{equation}\label{BLDdescription}
x\mapsto B(x) + \tilde{Y}(x).
\end{equation}
As before, $B$ denotes a Brownian motion started at zero and with variance $1/2$ and mean $0$ at time $1$. The process $\tilde{Y}$ is independent of $B$ and has a Radon-Nikodym derivative against Brownian motion started at zero and with variance $1/2$ and mean $0$ at time $1$ which is propositional to
\begin{equation}\label{BLDdescriptionRND}
\left(\int\limits_0^1 e^{-2\tilde{Y}(s)}ds\right)^{-u}\left(\int\limits_0^1 e^{2\tilde{Y}(1)-2\tilde{Y}(s)}ds\right)^{-v}.
\end{equation}
The factor $e^{-2v\tilde{Y}(1)}$ which appears above can be thought of as introducing a drift so we can also characterize $\tilde{Y}$ as having a Radon-Nikodym derivative against Brownian motion started at zero and with variance $1/2$ and mean $-v$ at time $1$ which is propositional to
\begin{equation}\label{BLDdescriptionRND2}
\left(\int\limits_0^1 e^{-2\tilde{Y}(s)}ds\right)^{-u-v}.
\end{equation}

It may not seem so immediate how to move between the two descriptions above for $h_{u,v}$. Such a matching is shown in \cite{bryc2021markov}. Here is a brief sketch. Starting with \eqref{BBKWdescription}, write $\hat{Y}(x) = Y(0)-Y(x)$. The aim is to show that $\hat{Y}(\cdot)$ and $\tilde{Y}(\cdot)$ have the same law. Notice that the Radon-Nikodym derivative from \eqref{BBKWdescriptionRND} can be rewritten in terms of $\hat{Y}$ and $Y(0)$ as
$$
\exp\left(2 (u+v) Y(0) - e^{2Y(0)} \int\limits_0^{1} e^{-2 \hat{Y}(s)}ds \right) e^{-2v \hat{Y}(1)}.
$$
The law of $\hat{Y}(\cdot)$ and $Y(0)$ are independent, and hence it is possible to integrate out $Y(0)$ since this does not figure into the description of $h_{u,v}$. Recall the integral identity
$$
\int_{-\infty}^{\infty} e^{2a x-b e^{2x}}dx = \tfrac{1}{2} b^{-a}\Gamma(a)
$$
which holds provided the real part of $a$ and $b$ are strictly positive. Choose $a=u+v$ and $b = \int_0^{1} e^{-2\hat{Y}(s)}ds$. Since both $a$ and $b$ are strictly positive, after integrating this yields the Radon-Nikodym derivative expression in \eqref{BLDdescriptionRND} (written there in terms of $\tilde{Y}$, not $\hat{Y}$) and hence implies that $\hat{Y}$ and $\tilde{Y}$ have the same law.

As far as how to link either of these descriptions to the Laplace transform formula, I will just say a bit since this is a substantial calculation. As is generally the case, given the description of the process $h_{u,v}$ as above, it is much easier to verify that it has the correct Laplace transform than it is to invert the Laplace transform from the start without any inspiration. Taking the Laplace transform relies on the following ideas (much more detail can be found in the papers \cite{BKWW} and \cite{BLD}). First, recall that for Brownian motions subject to energetic penalization by an exponential potential, the submarkov generator $\mathcal{L}$ that is relevant acts on a suitable space of functions $f$ as
$\mathcal{L}f:= f''(x) - e^{2x}f(x)$.
Define  $f_{u}(x):= K_{\iota u}(e^x)$ where $K$ is the modified Bessel function (or Macdonald function) given, for $u\in \mathbb{R}$ and $x>0$, by $K_{\iota u}(x) = \int_0^{\infty} e^{-x\cosh(w)}\cos(uw)dw$. Then the heat kernel from $x$ to $y$ in time $t$ for this operator takes the form
$$
\int\limits_0^{\infty} e^{-t u^2} f_u(x)f_u(y) \frac{2 du}{\pi |\Gamma(\iota u)|^2}
$$
as can be shown by appealing to the eigenrelations $\big(\mathcal{L} f_u\big)(x) = -u^2 f_u(x)$. This leads to formulas for the finite dimensional distributions of $h_{u,v}$ in terms of integrals of these functions. Finally, the Laplace transform can be computed by appealing repeatedly to the two identities:
\begin{align*}
\int\limits_{-\infty}^{\infty} e^{tx} f_u(x)f_v(x)dx
&= \frac{2^{t-3}}{\Gamma(t)} \left|\Gamma\left(\frac{t+\iota (u+v)}{2}\right)\Gamma\left(\frac{t+\iota (u-v)}{2}\right)\right|^2,\\
\int\limits_{-\infty}^{\infty} e^{tx} f_u(x)dx
&= 2^{t-2}\left|\Gamma\left(\frac{t+\iota u}{2}\right)\right|^2,
\end{align*}
where $u,v\in \mathbb{R}$ and $\textrm{Re}(t)>0$. The gamma function products arising here ultimately account for those in the Laplace transform formula.

\subsection{A few comments and directions}
A keen reader will have noticed that the Laplace transform formula in \cite{CK} is only given for $u+v>0$. I will discuss later about the origins of this restriction. The natural question here is to identify the Laplace transform or better yet the law of $h_{u,v}$ when $u+v<0$. In fact, it is even non-trivial to take the limit where $u+v\searrow 0$ in the Laplace transform formula \eqref{eq:laplace} or the formulas in  \eqref{BBKWdescription} and  \eqref{BBKWdescriptionRND} for $h_{u,v}$ from \cite{BKWW}. At least in the case of \eqref{eq:laplace}, this limit can be performed with some care and does recover the Brownian motion Laplace transform as one would hope. By comparison, the formula for $h_{u,v}$ in \eqref{BLDdescription} and \eqref{BLDdescriptionRND} (or even better \eqref{BLDdescriptionRND2}) immediately admits this limit -- when $u+v=0$ in \eqref{BLDdescriptionRND2} the Radon-Nikodym derivative is identically equal to $1$ and hence $\tilde{Y}$ is a Brownian motion of standard deviation $\sqrt{2}$ and drift $-4v$. Plugging this into \eqref{BLDdescription} shows that $h_{u,v}$ is a standard Brownian motion with drift $u=-v$, as it should be.

The description for $h_{u,v}$ from \cite{BLD} in \eqref{BLDdescription}, in fact, readily admits and extension to all $u$ and $v$, not just $u+v>0$. In that case, the process $\tilde{Y}(\cdot)$ is still well-defined. What is not clear (though is conjectured to be true in \cite{BLD}) is whether the process $h_{u,v}$ defined through \eqref{BLDdescription} in this case when $u+v<0$ is, in fact, stationary for the KPZ equation. Proving (or disproving) this seems like a great question. A natural approach might be to try to use some variant of uniqueness of analytic continuation -- to
show that the stationary measure depends in some sort of analytic manner on the boundary parameters $u$ and $v$, and likewise $h_{u,v}$ defined above by \eqref{BLDdescription}. This is all rather vague since the correct notion of analyticity as well as the means to prove it is currently unclear to me.

Another extension discussed in \cite{BLD}  is to consider the open KPZ equation stationary measure on a general interval $[0,L]$, or even on a half-infinite interval $[0,\infty)$. In both cases they write down candidates for the stationary measure. It should be possible to modify the approach from \cite{CK} to address these cases.  One nice observation in \cite{BLD} is that the stationary measure process $h_{u,v}$ actually arose (in a very different context) in 2004 work of Hariya and Yor \cite{HariyaYor}. In that case, their motivation was to study properties of exponential transforms of Brownian motion.

On the subject of the form of the Radon-Nikodym derivatives above, I just want to mention that they bare remarkable similarities to the type of reweighting that arose in my study with Hammond of the KPZ line ensemble \cite{CHammond}.

My final remark returns to the question of uniqueness (for each pair $u,v$) of the open KPZ stationary measures $h_{u,v}(\cdot)$. In \cite{CK} we conjectured that this is true, and moreover that the open KPZ equation is ergodic and satisfies a one-force one-solution principle which essentially says that if you start at time $-T$ with two different choices of initial data, and then look at the solution around time $0$, as $T\to \infty$, the height increments of the two processes will converge to be the same, and hence independent of the initial data. There are some similar results to this proved in the literature for the periodic boundary condition KPZ equation or for some other similar models, see for example \cite{EKMS, BCKFullLine, hairer2018, alex2019stationary, gubinelli2018infinitesimal, rosati2019synchronization}.

\section{Open ASEP to KPZ}\label{sec:aseptokpz}

How does one even start to show that there exists a stationary measure for the open KPZ equation, let alone check that it is given in the form that we claim? Usually the idea is to construct the generator or semi-group for the open KPZ equation and then check that the generator acts on a given stationary measure to give zero or that the semi-group preserves the stationary measure. Mixing properties can then be studied by finding the spectral gap or probing clever couplings.

This approach has been successfully implemented for the KPZ equation on a torus or to some extent on the full line where the stationary measure is purely Brownian, see for example \cite{FQ14,hairer2018,gubinelli2018infinitesimal}. For the open KPZ equation this approach has not been implemented, though I would be quite interested to see it done.

Instead, in \cite{CK} we proceed through a discretization of the open KPZ equation -- the interacting particle system called {\it open ASEP}. In a nutshell, the idea is to first show that under special scaling, the height function for open ASEP converges to the open KPZ equation, provided that the initial data has a limit and satisfies some reasonable hypotheses (i.e., its has Brownian-like H\"older behavior). Open ASEP is a finite state space Markov process and for each system size $N$ it has a unique stationary measure. The challenge then becomes to show that these stationary measures converge to a limit as $N$ goes to infinity and satisfy the desired hypotheses. In fact, we only show the existence of a limit in the case where $u+v>0$ (I will explain what these mean in the open ASEP context below). For  general $u,v$ we are able to show tightness of the $N$-indexed sequence of open ASEP stationary measures which translates into the existence of subsequential limits that are stationary measures for the open KPZ equation. Of course, if we knew uniqueness of the open KPZ stationary measures, that would imply there is only one limit point.

In this section, I am going to try to explain how open ASEP approximates the open KPZ equation. This is based on work of mine with Shen in \cite{CS} and a subsequent extension by Parekh \cite{Parekh}. Open ASEP is quite an interesting and well-studied object in its own right, and the rest of this note will almost entirely focus on it. As such, I will start out here with a bit of background, in particular its phase diagram. Then I will return to the connection to open KPZ.

\subsection{Introducing open ASEP}
Open ASEP was first introduced by MacDonald, Gibbs and Pipkin in 1968 \cite{MGP} as a model for  the dynamics of ribosomes on an mRNA chain during the synthesis of proteins. Already in that work, their main interest was in studying its stationary measure. Within probability, Spitzer initiated study of a general class of exclusion processes in his 1970 work \cite{SPITZER1970246} and Liggett introduced open ASEP to the community in his 1975 work \cite{Liggett75} as a tool to study the nature of stationary measures for ASEP on the full line and half-line.

\begin{figure}[h]
\centering
\scalebox{0.8}{\includegraphics{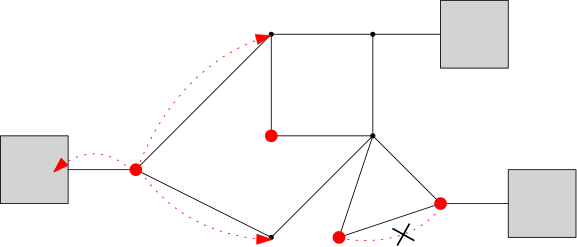}}
 \captionsetup{width=.9\linewidth}
\caption{An exclusion process on a network with reservoirs. Dotted red lines show the possible jumps of the left-most particle and the dotted red line with an $\times$ through it, shows an excluded jump. Other possible jumps and jump rates are not shown.}
\label{Fig:GeneralExclusionProcess}
\end{figure}

Let me now introduce open ASEP. Consider the exclusion process in Figure \ref{Fig:GeneralExclusionProcess}. Particles (red dots) occupy vertices of a graph and jump along edges according to exponential rates subject to the exclusion rule (jumps to occupied sites are suppressed). There are {\em reservoirs} (grey squares) which (according to exponential rates) insert particles into unoccupied neighboring sites, or remove particles from occupied neighboring sites. This models transport through a network with sources and sinks.

The presence of reservoirs leads to remarkable physical behavior. Since the total number of particles is not conserved, even when in its steady state (i.e., started its unique stationary measure) there will typically be a net flow  of particles through this system---like how a flowing stream may have a stationary density profile while still moving water from source to sink. In statistical physics such models are said to have a {\it non-equilibrium steady state} (e.g. see the review \cite{Blythe_2007}). Because the Markovian dynamics which produce these steady states are not reversible, non-equilibrium steady states do not take the Boltzmann weight form common in equilibrium statistical mechanics.  This renders the study of thermodynamic (i.e. large system size) limits of non-equilibrium steady states quite challenging.

The presence of reservoirs can also induce phase transitions as the number of nodes ($N$ in our case) in the network grows. I am unaware of a precise statements to this effect when dealing with general networks but for the one-dimensional asymmetric simple exclusion process with reservoirs (i.e., open ASEP) such a phase transition is understood. Besides serving as a transport model, I will explain below how the height function for open ASEP is also connected with stochastic interface growth and, through exponentiation, to a discrete stochastic heat equation.

\begin{figure}[h]
\centering
\scalebox{0.8}{\includegraphics{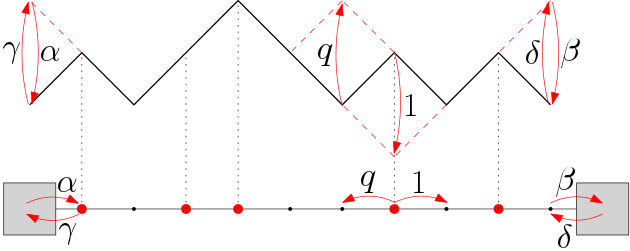}}
 \captionsetup{width=.9\linewidth}
\caption{A given instance of open ASEP with $N=10$ is illustrated both in terms of occupation variables (the red dots) and its height function $h_N$ (the piecewise linear function drawn above the particle configuration). The red arrows indicate some possible moves with the associated exponential rates labeled.}
\label{Fig:openASEP}
\end{figure}

Let me start by defining the open ASEP in one-dimension. Figure \ref{Fig:openASEP} gives an illustration of it along with its height function. We consider an $N$ site nearest neighbor graph $\{1,\ldots, N\}$ with edges between consecutive numbers. Particles jump left according to exponential clocks of rate $q<1$ and right at rate $1$. These jumps (and all others) are independent and only taken if the exclusion rule is observed -- particles cannot move to occupied sites. In addition to this bulk evolution, there are two reservoirs which interact with site $1$ and site $N$. At rate $\alpha$ a particle enters site $1$ (provided it is empty) and at rate $\gamma$ a particle is removed from site $1$ (provided it is occupied); similarly at rate $\delta$ a particle enters site $N$  (provided it is empty) and at rate $\beta$ a particle is removed from site $N$ (provided it is occupied).
This process can be encode as $\tau(t)=\big(\tau_1(t),\ldots, \tau_N(t)\big)$ where $\tau_x(t)=1$ if site $x$ is occupied at time $t$ and otherwise $\tau_x(t)=0$ if it is unoccupied. The (backward) generator of the process $\tau(t)$ is discussed later around \eqref{eq:genzero}. It is important to note that since there is no conservation of particle number (particles are created and removed) there is a single invariant measure (or steady state) for this process.
The physics literature often denotes this measure as $\pi_{N}(\tau)$, with the dependence on the other parameters implicit, and denotes the expectation of a function $f:\{0,1\}^{\llbracket 1,N\rrbracket}\to \R$ under $\pi_{N}$ by
\begin{equation}\label{eq:lrangle}
\big\langle f\big\rangle_N:= \!\!\!\!\sum_{\tau\in \{0,1\}^{\llbracket 1,N\rrbracket}} f(\tau)\cdot \pi_N(\tau).
\end{equation}
While the existence of this measure $\pi_{N}$ is clear, it is not obvious how it behaves as $N$ tend to infinity. I will return to this important point a bit later.

There is one piece of information that is  lost in the open ASEP occupation process -- the count of how many particles have entered or exited from the boundaries. Going to the height function process remedies this. I will define the  height function process $h_N(t,x)$ as follows. Of course, the subscript indicates the lattice size $N$, and the time $t$ and spatial location $x$ are now both written as arguments. The dependence of $h_N(t,x)$ on the other parameters $q,\alpha,\beta,\gamma,\delta$ will be generally suppressed from the notation, though they will eventually all depend on $N$ non-trivially. The height function is defined for $t\geq 0$ and $x\in\{0,\ldots N\}$ as
\begin{equation}\label{eq_hN}
h_N(t,x) := h_N(t,0) + \sum_{i=1}^{x} \big(2\tau_i(t)-1\big),\qquad \textrm{where}\quad h_N(t,0):=-2\mathcal{N}_N(t)
\end{equation}
and where the net current $\mathcal{N}_N(t)$ records the total number of particles that have entered into site $1$ from the left reservoir up to time $t$ minus the number of particles that have exited from site $1$ into the left reservoir up to time $t$. It is convenient to linearly interpolate to define a continuous height function in space.

The open ASEP height function process does not have an invariant probability measure. Indeed, due to the net current $\mathcal{N}_N(t)$, the only invariant measure takes the form of infinite counting measure on the height at the origin and then the induced measure (coming from $\pi_{N}$ on the occupation variables $\tau$) on the height function increments from there. Another way of saying this is that if $h_N(x)$ is the random height function defined by $h_N(x):=\sum_{i=1}^{x} \big(2\tau_i(t)-1\big)$ where $\tau$ is distributed according to the invariant measure $\pi_N(\tau)$, then starting the open ASEP height process with initial data $h_N(0,\cdot) =h_N(\cdot)$ implies that for all later times $t>0$, $h_N(t,x)-h_N(t,0)$ will still have the law of $h_N(x)$ as a process in $x$. Thus, I call $h_N(\cdot)$ the stationary measure. Here I use stationary instead of invariant to emphasize that it is the increment process that is invariant.

\subsection{Phase diagram}\label{sec:phasediagram}
The boundary reservoirs play a key role in determining the limiting behavior of open ASEP as $N$ goes to infinity.
Define the current in stationarity to be
$$
J_N:= \frac{\big\langle \alpha(1-\tau_1)-\gamma\tau_1 \big\rangle_N}{1-q}.
$$
The term $\alpha(1-\tau_1)$ accounts for the rate $\alpha$ at which particles enter the system provided $\tau_1=0$ and the term $\gamma\tau_1$ for the rate $\gamma$ at which they depart provided $\tau_1=1$. The different measures the instantaneous rate of signed movement across the bond between the reservoir and site $1$.

If open ASEP  starts according to its invariant measure $\pi_N$, then it follows that for all $t>0$
$$
J_N= \frac{\big\langle \mathcal{N}_N(t)\big\rangle_N}{t(1-q)}.
$$
I am abusing notation here since now $\langle\cdot\rangle_N$ represent the expectation of the open ASEP occupation process $\tau(t)$ started from in its invariant measure $\pi_N(\cdot)$ at time $t=0$, not just the expectation of a function of the state space.
The quantity $J_N$ represents the average (normalized by the bulk drift $1-q$) current of particles moving through the system at stationarity. In fact, since particles are conserved by the bulk dynamics, this average current will be the same everywhere in the system.

There is a remarkable phase diagram for the large $N$ limit of $J_N$ which highlights the role of the boundary rates. The existence of such a limit is non-trivial, let alone computing its value as a function of $q,\alpha,\beta,\gamma,\delta$. I will give a simple (and non-rigorous) heuristic to derive it, though a more complete derivation in the physics literature follows from the matrix product ansatz (which will be introduced near the end of this note) in \cite{Derrida_1993} for $q=0$ and \cite{USW} for general $q$.
In fact, a version of this phase diagram arose in a closely related context in early work of \cite{Liggett77}.

The limit $J:=\lim_{N\to \infty} J_N$ exists and depends on two parameters $\rho_{\ell}$ and $\rho_r$ which have nice physical interpretations as effective boundary densities at the left and right boundaries.


Consider $N$ very large and focus on the invariant measure near site $1$. Without any justification, imagine that locally and asymptotically in $N$, the invariant measure looks like a product Bernoulli measure there with density $\rho_\ell$ for some $\rho_\ell\in (0,1)$. Provided this, it is possible to determine what value $\rho_\ell$ must take through a simple consideration. Since particles are conserved within the bulk of open ASEP, in stationarity the net number of particles moving from the left reservoir into site $1$ must equal the net number moving from site $1$ to $2$. Taking expectations, this implies a simple conservation equation
$$
\alpha (1-\rho_\ell) -\gamma \rho_\ell = (1-q) \rho_\ell(1-\rho_\ell).
$$
The first term on the left accounts for the rate $\alpha$ at which particles enter from the boundary to site $1$ provided it is unoccupied (which happens in stationarity with probability $1-\rho_\ell$) and the second term accounts for the rate $\gamma$ at which particles exit to the boundary from site $1$ provided it is occupied (which happens in stationarity with probability $\rho_\ell$). On the right, the factor $\rho_\ell(1-\rho_\ell)$ is the probability (under the Bernoulli product measure assumption) of having a particle and a hole next to each other in any prescribed order. The factor $(1-q)$ comes from the rate of particles jumping right minus the rate jumping left. Solving this quadratic equation yields $\rho_\ell$. A similar consideration around $N$ yields $\rho_r$.  In both cases, there is just one positive solution which yields the effective density at the boundaries.

It is convenient (for later purposes) to write the solutions to the above quadratic equations explicitly in the following form:
\begin{equation}
\rho_\ell = \frac{1}{1+C},\qquad \rho_r = \frac{A}{1+A}
\end{equation}
where $A, C>0$ are given by
\begin{equation}\label{eq_ABCD}
A=\kappa^+(q,\beta,\delta), \quad B=\kappa^-(q,\beta,\delta),\quad C=\kappa^+(q,\alpha,\gamma) ,\quad D=\kappa^-(q,\alpha,\gamma)
\end{equation}
with
\begin{equation}\label{eqnkappa}
\kappa^{\pm}(q,x,y):=\dfrac{1}{2 x}\left(1-q-x+y\pm\sqrt{(1-q-x+y)^2+4 xy} \right).
\end{equation}
The $\kappa^{-}$ terms are the negative roots of the quadratic and are introduced for later purposes. Notice that for $q$ fixed, \eqref{eq_ABCD} actually gives a bijection between $\{(\alpha,\beta,\gamma,\delta):\alpha,\beta>0,\gamma,\delta\geq 0\}$ and $\{(A,B,C,D):A,C>0,B,D\in (-1,0]\}$.

\begin{figure}[h]
\centering
\scalebox{0.4}{\includegraphics{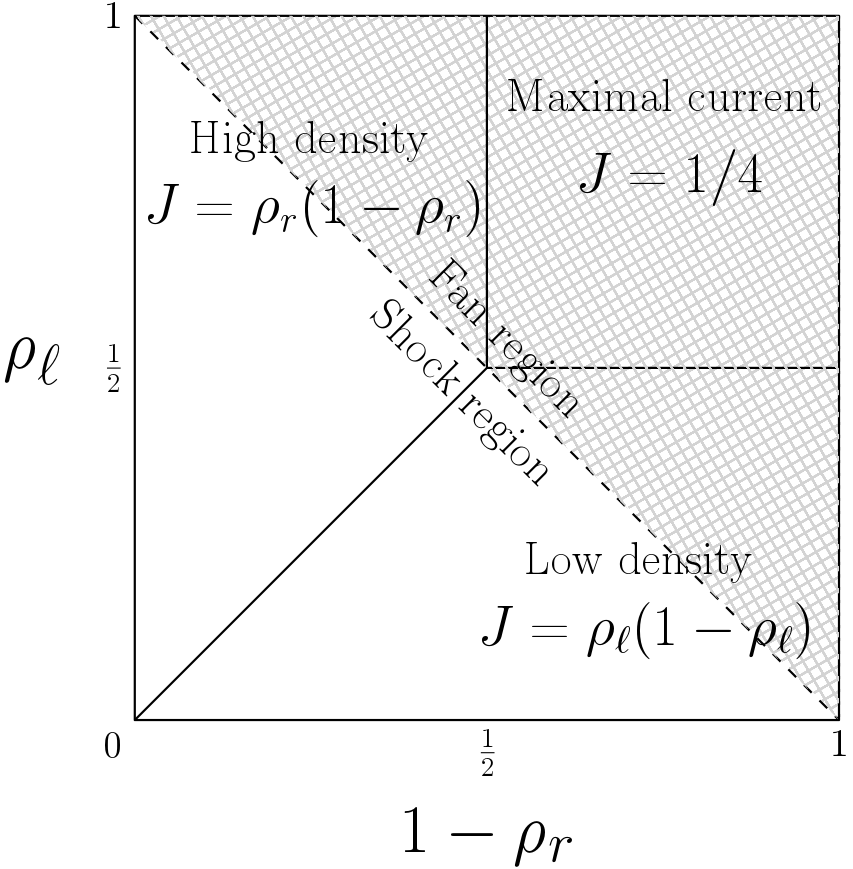}}
 \captionsetup{width=.9\linewidth}
\caption{The phase diagram for the current $J$ of open ASEP as a function of the effective densities $\rho_\ell$ and $\rho_r$.}
\label{Fig:phasediagram}
\end{figure}

The phase diagram for the value of $J$ is determined entirely by the values of $\rho_\ell$ and $\rho_r$ (or equivalently $A$ and $C$) as follows (see also Figure \ref{Fig:phasediagram}):
\begin{itemize}[leftmargin=*]
\item {\em Maximal current phase ($\rho_\ell>1/2$ and $1-\rho_r>1/2$, or equivalently $C<1$ and $A<1$)}: The left boundary creates particles at a fast enough rate (and does not remove them too quickly) and the right boundary removes particles at a fast enough rate (and does not create them too quickly) so that the system is able to transport particles from left to right in the bulk at its level of maximal efficiency. The maximal rate of transport for particles is $1/4$ and is achieved when the bulk density is $1/2$. To see this, note that if there is local product measure of density $\rho$ in the bulk of ASEP then the rate at which particles move (normalized by $1-q$) will be $\rho(1-\rho)$. This is maximized at $\rho=1/2$ and takes value $1/4$ in that case. Thus, $J=1/4$ in this phase and the density should look locally in the bulk like product Bernoulli with density $1/2$.
\item {\em Low density phase ($\rho_\ell+\rho_r<1$,  or equivalently $C<A$ and $A>1$)}: The left boundary creates particles relatively slowly and the right boundary removes them fast enough so that they do not build up there. As such, the density in the bulk is determined by the effective density on the left, $\rho_\ell$ and the current $J=\rho_\ell(1-\rho_\ell)$.
\item {\em High density phase ($\rho_\ell+\rho_r>1$, or equivalently $C>A$ and $C>1$)}: By applying a particle-hole transform this is equivalent to the low density phase for the holes. In particular, now the right boundary removes particles so slowly and there is enough input from the left boundary that the density builds up and becomes $\rho_r$. Thus, the current $J=\rho_r(1-\rho_r)$.
\end{itemize}
The line $\rho_\ell+\rho_r=1$ involves coexistence of the high and low density phases and displays interesting behavior that I will not touch on here. Another point of interest is the {\em triple point} when $\rho_\ell=\rho_r=1/2$. This point will play an key role since it is around there that the KPZ equation arises.

Besides the phases addressed above, there is one other important division of the phase diagram:
\begin{itemize}[leftmargin=*]
\item {\em Fan region ($\rho_\ell>\rho_r$, or equivalently $AC<1$)}: Very close (in a scaling going to zero relative to $N$) to the left boundary, the density exceeds that of the bulk. This is called the fan region since going from high to low density in the Burgers equation produces a rarefaction fan.
\item {\em Shock region ($\rho_\ell<\rho_r$, or equivalently $AC>1$)}: Very close (in a scaling going to zero relative to $N$) to the left boundary, the density is lower than that of the bulk. Thus, in the same spirit of the Burgers equation, one sees shock-type behavior here.
\end{itemize}

The boundary between the fan and shock regions, when $\rho_\ell=\rho_r$ or equivalently $AC=1$, is special since it is the only part of the phase diagram where the invariant measure is simple. Along this line, in the finite $N$ open ASEP, the invariant measure is Bernoulli product with density $\rho=\rho_\ell=\rho_r$. I will recall later how this important fact is shown.

Beyond determining the current $J$, the above phase diagram also dictates the nature of the fluctuations and large deviations for the open ASEP height function stationary measure and process. For instance, the fluctuations as a spatial process for the stationary measure has been determined for TASEP ($q=0$) by Derrida, Lebowitz and Enaud \cite{DEL} and for ASEP (general $q\in (0,1)$) more recently by Bryc and Wang \cite{BW}. In fact, the work of \cite{BW} is what prompted my interest in the open KPZ stationary measure and what indicated to me that it should be possible using, in part, methods in that work.

\subsection{Microscopic Hopf-Cole transform}

I want to now address the question of how open ASEP and open KPZ are related. The study of stochastic PDE limits of interacting particle systems has been a hot topic in the past decade, and has a significantly longer history. Typically such limits arise when the system size and time scale are taken to infinity appropriately while various parameters are tuned in a critical manner as well. In the case of the KPZ equation on the entire real line, the earliest example of such a convergence result is the 1995 work of Bertini and Giacomin \cite{cmp/1158328658}. That work relied on the starting observation which came from earlier work of G\"{a}rtner \cite{Gartner} that
\begin{equation}\label{eq:MHC}
Z_t(x) := e^{-\lambda h_N(t,x) + \nu t},\qquad \textrm{with} \quad \lambda = \tfrac{1}{2}\log q,\quad \nu = 1+q+-2\sqrt{q}
\end{equation}
satisfies a microscopic version of the SHE. This transform can be thought of as a microscopic version of the Hopf-Cole transform. The analysis from there on was challenging, but possible since a well-developed solution theory for the SHE was already developed.

I have a few papers from the past decade devoted to generalizing the Bertini-Giacomin approach in different directions and to different types of interacting particle systems. Always, the starting point is the microscopic Hopf-Cole transform. Often this arises as a result of a hidden Markov duality, though I will not say more about this here. The need for such a transform severely limits the applicability of this technique. However, work of Dembo and Tsai \cite{Dembo2016}, and then subsequent developments by Yang have demonstrated that even in instances where the microscopic Hopf-Cole transform does not exactly satisfy a version of the SHE, it is still possible to show convergence to the SHE, and hence of the height process to the KPZ equation. In fact, Yang's recent work \cite{yang2021kardarparisizhang} shows that the open KPZ equation arises as a limit of a wide variety of open exclusion processes. Other techniques such as energy solutions, paracontrolled distributions and regularity structures also provide routes to derive the KPZ equation as a limit of interacting particle systems.

In the case of the open ASEP, the first proof of its convergence to the open KPZ equation came in my work with Shen  \cite{CS}. As indicated above, the starting point was our observation that open ASEP satisfies an exact microscopic Hopf-Cole transform. In fact, a special case of this observation showed up around the same time in work of Gon\c{c}alves, Landim and Milan\'{e}s \cite{10.1214/16-AAP1200}.

In order for the microscopic Hopf-Cole transform to work for open ASEP, it is necessary to impose a restriction on the parameters. With Shen, we called this the {\it Liggett condition} since it arose in his much earlier work \cite{Liggett75,Liggett77}. The condition is that
\begin{equation}\label{eq:Liggett}
\frac{\alpha}{1} + \frac{\gamma}{q} = 1 = \frac{\beta}{1}+\frac{\delta}{q}.
\end{equation}
Under this condition the effective reservoir densities simplify so that $\rho_\ell = \alpha$ and $\rho_r = 1-\beta$. Moreover, this condition implies that the role of the reservoirs can be replaced by a simpler boundary interaction -- if there is no particle at site 1, then with rate $1$
the system attempts to pull a particle out of the reservoir to fill site 1, and that attempt is successful with probability $\alpha$; if there is a particle at site $1$, then with rate $q$ the system attempts to push that particle out into the reservoir, and that attempt is successful with probability $1-\alpha$. There is a similar dynamic with the reservoir near $N$. Thus, under Liggett's condition, the left reservoir can be thought of as a standard ASEP site which is occupied or vacant with probability $\alpha$ or $1-\alpha$ at any given moment of time, independently of all other times. This implication of the condition was explained initially in \cite{Liggett75}.

Under Liggett's condition, $Z(t,x)$ from \eqref{eq:MHC} satisfies the following discrete stochastic differential equation version of the SHE:
\begin{equation}\label{eq:dSHE}
d Z_t(x) = \tfrac{1}{2} \Delta Z_t(x) + dM_t(x)
\end{equation}
for all $x\in \{0,\ldots, N\}$ subject to the boundary conditions that for all $t\geq 0$
$$
Z_t(-1) = \mu_\ell Z_t(0)\qquad Z_t(N+1) = \mu_r Z_t(N)
$$
where $\mu_\ell$ and $\mu_r$ take values in $[q^{1/2},q^{-1/2}]$ and are given by
$$
\mu_\ell = q^{-1/2} - \alpha(q^{-1/2}-q^{1/2})\qquad\textrm{and}\qquad \mu_r = q^{-1/2} - \beta(q^{-1/2}-q^{1/2}).
$$
The term $M_t(x)$ represents a martingale with an explicit bracket process that I will not record here. The meaning of the boundary condition (which is an inhomogeneous discrete Robin boundary condition) is that when considering $\Delta Z_t(x)$ for $x=0$ or $x=N$, the term which involves $x=-1$ or $x=N+1$ is replaced by use of the boundary condition above.

The convergence result that Shen and I proved, and that will be described in the next section, uses this transformation as the starting point and thus requires that the parameters satisfy Liggett's condition. This reduces the number of boundary parameters from four to two, though these remaining two parameters can be thought of as the reservoir potentials $\rho_{\ell}$ and $\rho_r$ and by tuning them we are able to access the full two-parameter family of open SHE/KPZ equation boundary parameters. It would be nice to see a proof of the open ASEP to KPZ scaling limit which does not require Liggett's condition (note that this condition seems to be necessary in Yang's work \cite{yang2021kardarparisizhang} as well).

\subsection{Convergence to open KPZ}\label{sec:conv}

In order to prove convergence of the discrete SHE \eqref{eq:dSHE} to the continuum SHE, it is necessary to introduce two more scalings of the remaining parameters (which we can take to be $q,\rho_\ell$ and $\rho_r$). The first is known as weak asymmetry scaling. In the full-line setting of Bertini and Giacomon, they introduce a scaling parameter $\epsilon$ and take $q$ scaled close to $1$ on the order of $\epsilon^{1/2}$, space scale down by $\epsilon$ and time scaled down by $\epsilon^{2}$. It takes a bit of work to see why this is a natural choice of scaling and, in my opinion, is best seen by studying scalings under which the KPZ equation remains fixed (e.g., see \cite{ICAMS}).

In the context of the open ASEP, space is always fixed to involve $N$ sites. This suggests to think of $N$ and $\epsilon^{-1}$ as being the same (or at least proportional). Informed by this and earlier scaling of Bertini and Giacomin, I will assume below that
$$
q=\exp\left(-\frac{2}{\sqrt{N}}\right)
$$
and call this weak asymmetry scaling. Notice that $q\approx 1-2N^{-1/2}$, which matches with the $\epsilon^{1/2}$ notion of weak asymmetry coming from Bertini and Giacomin. Likewise, further define
$$
h^{(N)}(t,x):= N^{-1/2} h_N(\tfrac{1}{2} e^{N^{-1/2}}N^2 t, Nx) + (\tfrac{1}{2}N^{-1} + \tfrac{1}{24})t,\qquad z^{(N)}(t,x) := e^{h^{(N)}(t,x)}.
$$
The scaling of space like $Nx$ means that in $h^{(N)}(t,x)$ and $ z^{(N)}(t,x)$, $x$ can vary over $[0,1]$. The scaling of time like $\tfrac{1}{2} e^{N^{-1/2}}N^2 t$ captures the $N^2$ scaling from Bertini and Giacomin, as well as some additional corrections that simplify coefficients elsewhere. The rest of the definition is fixed by the microscopic Hopf-Cole transform -- namely the fact that $z^{(N)}(t,x)$ should satisfy a scaled version of the discrete SHE. Indeed, the factor $ (\tfrac{1}{2}N^{-1} + \tfrac{1}{24})t$ comes from the expansion of $\nu \tfrac{1}{2} e^{N^{-1/2}}N^2 t$ down to order $1$ terms in $N$.

So far, none of this scaling has involved the boundary parameters or conditions. Since $h^{(N)}(t,x)$ involves diffusive scaling between space and the scale of the height function it is natural to imagine that the density of particles in open ASEP should be close to $1/2$ in order for this scaling to make sense. In fact, it should be within order $N^{-1/2}$ of density $1/2$ to respect this scaling. This implies that $\rho_\ell$ and $\rho_r$ should be scaled in this manner around density $1/2$, i.e., in an $N^{-1/2}$ window around the triple point of the phase diagram. Specifically, assume now that for some $u,v\in \mathbb{R}$,
$$
\rho_\ell = \frac{1}{2} + \frac{u}{2} N^{-1/2} + o(N^{-1/2}),\qquad \rho_r = \frac{1}{2} - \frac{v}{2} N^{-1/2} + o(N^{-1/2}).
$$

The above assumptions have all been on parameters, but it is also necessary to assume something about the initial data. Since $N$ is varying, for each $N$ there will be a different choice of initial data -- denote the height function for that initial data by $h_N(\cdot)$. This is an abuse of notation from earlier where I let this denote the stationary height function. For the moment, I will just take this to be any choice of initial data. Introduce its scaled version $h^{(N)}(x):=N^{-1/2} h_{N}(Nx)$ and exponential $Z^{(N)}(x):=e^{h^{(N)}(x)}$. These definitions match with the $t=0$ height function scaling introduced above.

The following {\it H\"older bounds} assumption will be sufficient to state the KPZ convergence result: For all $\theta\in (0,1/2)$ and every $n\in \mathbb{Z}_{\geq 1}$ there exists $C(n),C(\theta,n)>0$ such that for every $x,x'\in [0,1]$ and $N\in \mathbb{Z}_{\geq 1}$
$$
||Z^{(N)}(x)||_{n} \leq C(n),\qquad \textrm{and}\qquad ||Z^{(N)}(x)-Z^{(N)}(x')||_n \leq C(\theta,n)|x-x'|^{\theta}
$$
where $||\cdot||_n := \mathbb{E}[|\cdot|^n]^{1/n}$ and $\mathbb{E}$ is the expectation over the initial data.

Notice that assuming initial data satisfying these H\"older bounds, the initial data $Z^{(N)}(\cdot)$ will form a tight sequence of random functions. With some work (as done in \cite{CS}) it can be shown that the discrete SHE preserves this class of initial data. This is the first step to proving convergence to the continuum SHE/KPZ equation since it implies tightness of the entire process. The second step is to show that all subsequential limits satisfy the desired SHE/KPZ equation.

I will informally state the main convergence result from \cite{CS} and \cite{Parekh}. This was proved in \cite{CS} under the assumption that $u,v\geq 1/2$ and extended by \cite{Parekh} to general $u,v\in \mathbb{R}$ case. The reason why the general case was harder is that the boundary condition for $u,v\geq 1/2$ can be thought of as repulsive or killing in terms of the Feynman-Kac representation, and thus the heat kernel for the associated Laplacian tends to decay to zero, simplifying various estimates. In the general case, the heat kernel may grow exponentially, thus complicating matters.

The combination of these two works showed that, provided a sequence of $N$-indexed open ASEP processes with parameters satisfying Liggett's condition, weak asymmetry and triple point scaling, and initial data satisfying H\"older bounds, then the following holds. For any fixed time horizon $T>0$, the law of $\{Z^{(N)}(\cdot,\cdot)\}_N$ is tight in the Skorokhod space $D([0,T],C([0,1]))$ of time-space processes that are CADLAG in time and continuous in space. Moreover, any limit point is in $C([0,T],C([0,1]))$, i.e., continuous in both time and space. If  there exists a (possibly random) non-negative-valued function $z_0\in C([0,1])$ such that $Z^{(N)}(\cdot)$ converges to $z_0(\cdot)$ along a subsequence as $N$ goes to infinity in the space of continuous processes on $[0,1]$, then along that same subsequence $Z^{(N)}(\cdot,\cdot)$ converges to $z(\cdot,\cdot)$ in $D([0,T],C([0,1]))$ where $z(\cdot,\cdot)$ is the unique (mild) solution to the SHE with boundary parameters $u$ and $v$ and initial data $z_0(\cdot)$, see Section \ref{sec:defopenKPZ}. In other words, this implies that the corresponding height function $h^{(N)}(\cdot,\cdot)$ converges to $h(\cdot,\cdot)$, the Hopf-Cole solution to the KPZ equation.

Before ending this section, I want to just remark on how the above scaling limit result relates to the mixing time conjecture for open ASEP and to the uniqueness conjecture for the stationary measure of open KPZ. Besides controlling the current and stationary measure density, the phase diagram for the open ASEP is supposed to control the mixing time, i.e., the time that it takes for a general initial state to converge (e.g. in total variation distance) close to the stationary measure. Recently, Gantert, Nestoridi and Schmid \cite{gantert2020mixing} have made progress on characterizing the mixing time. In the high and low density states, and assuming a strict asymmetry (i.e. $q<1$) they show that the mixing time grows linearly in the system size $N$. In the case of the triple point between all three phases, they are able to give an upper bound of order $N^3$. However, the expectation is that the true mixing time at this point and in the maximal current phase is of order $N^{3/2}$ and Schmid \cite{Schmid} has since proved this in the case of TASEP ($q=0$).

How does this ASEP mixing time behavior relate to the scaling limit of open ASEP to open KPZ? The mixing time of $N^{3/2}$ should have a prefactor that depends on the strength of the asymmetry. In particular, I expect it to behave like $(1-q)^{-1}N^{3/2}$. Under our weak asymmetry, this behaves like $N^{2}$ which is exactly the time scaling in order to arrive at the open KPZ equation. Thus, if someone can prove mixing for open ASEP in a manner sufficiently uniform over $q$, this could yield a route to study the mixing of the open KPZ equation as well (and hence imply the uniqueness of the stationary measure). Of course, this is probably not the easiest (or most direct) route to prove the uniqueness result about open KPZ.

\section{Taking asymptotics of the open ASEP stationary measure}\label{sec:asym}

The convergence result described above shows that under special scaling and assuming H\"older bounds on the open ASEP initial data, there is tightness of the scaled open ASEP height function and all subsequential limits of the initial data yield subsequential limits of the process which solve the open KPZ equation (in fact, things are phrased in terms of the SHE). Since my goal (i.e., the results claimed in Section \ref{sec:const}) is to construct a stationary measure for the open KPZ equation, I will consider now what happens when the convergence result is applied to the open ASEP stationary measure.

There are two types of results in Section \ref{sec:const} -- those that deal with the existence and general properties of the stationary measure and those that address the exact formulas in the case where $u+v>0$. Let me start by addressing the first type of result.

There are two inputs that we appeal to about the open ASEP stationary measure. The first is that when $\rho_\ell=\rho_r=\rho$, the corresponding stationary measure is Bernoulli with density parameter $\rho$. Liggett \cite{Liggett75} provided a proof of this (under the Liggett condition) and a more general result from the matrix product ansatz \cite{Derrida_1993}, as explained further in Section \ref{sec:DEHPrep}. It is straight-forward to see that in this case the open ASEP stationary measure height function satisfies the necessary H\"older bounds since it is just a simple random walk trajectory.

The second input is a microscopic version of the stochastic sandwiching result that I described earlier for the open KPZ stationary measure. Before stating it, let me motivate it. Imagine you have two versions of open ASEP, one with boundary parameters $\alpha,\beta,\gamma,\delta$ and the other with $\alpha',\beta',\gamma',\delta'$. If $\alpha\leq \alpha'$, $\beta\geq \beta'$, $\gamma\geq \gamma'$, and $\delta\leq \delta'$, then in the primed system there is an increased rate at which particles enter and a decreased rate at which they exit the system. It would reason then that the stationary measure for the primed system should have more particles than in the original system. This is true, as is the stronger statement that $\pi_N$ is stochastically dominated (see the definition below) by $\pi'_N$.

Consider two measures $\pi$ and $\pi'$ on $\{0,1\}^{\{1,\ldots, N\}}$. The measure $\pi$ is said to be stochastically dominated by $\pi'$ (written $\pi\preceq \pi'$) if there exists a coupling of $\pi$ and $\pi'$ on which all sites occupied under $\pi$ are likewise occupied under $\pi'$. More explicitly, this means that there exists a probability measure $\mu$ on $\{0,1\}^{\{1,\ldots, N\}}\times \{0,1\}^{\{1,\ldots, N\}}$ such that if let $(\tau,\tau')$ be sampled according to $\mu$ (here  $\tau$ and $\tau'$ take values in $\{0,1\}^{\{1,\ldots, N\}}$) then marginally $\tau$ has law $\pi$, $\tau'$ has law $\pi'$ and almost surely $\tau\leq \tau'$ in the sense that $\tau_i\leq \tau'_i$ for all $i\in \{1,\ldots, N\}$.

This result can be generalized to consider three sets of boundary parameters, and it is this consideration that leads to stochastic sandwiching for the open ASEP. After all, the height function increments are just sums of occupation variables, and thus stochastic domination of the individual occupation variables certainly implies that of their sums.

I will briefly explain how to demonstrate the above stochastic domination of the stationary measures. However, before that let me note how it implies the general $u,v$ H\"older bounds. When $u+v=0$ it is possible to choose $\rho_\ell=\rho_r$ and thus deal with product Bernoulli measure. When $u+v>0$ or $u+v<0$, the parameters $\alpha,\beta,\gamma,\delta$ can be adjusted in the spirit explained above to show that the corresponding stationary measure is stochastically sandwiched between two different product Bernoulli measures. The density parameter for the upper and lower bounding Bernoulli measures only differ by order $N^{-1/2}$ which is compatible with the diffusive scaling that is applied to the height function. Consequently, the H\"older bounds follow by applying the sandwiching in concert with the analogous bounds for simple random walks that are converging to drifted Brownian motions.

As far as proving the stochastic domination, the main idea is to use second class particles. This is the same mechanism that is used to show what is termed attractiveness of the usual ASEP.  Consider starting both the original and the primed version of open ASEP in state $\tau(0)$ and $\tau'(0)$ in such a way that $\tau(0)\leq \tau'(0)$ (e.g. they could both start out entirely empty). Let $\alpha(0) = \tau'(0)-\tau(0)$ denote the occupation variables for what are called {\it second class} particles. For ASEP on the full line, the basic coupling provides a dynamic on the pair $(\tau,\alpha)$ such that at any later time $\tau'(t)$ and $\tau(t)+\alpha(t)$ have the same distribution.

\begin{figure}[t]
	\captionsetup{width=.8\linewidth}
	\begin{center}
	\includegraphics[width=2in]{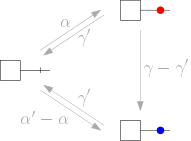}
	\end{center}
	\caption{Coupling $\tau$ and $\tau'$ at the boundary. Particles in $\tau$ are red dots, and those in $\tau'$ are blue dots. The transition rates and transitions are labeled in grey.}
	\label{fig:coupling}
\end{figure}

In order to demonstrate such a coupling on the interval, the basic coupling needs to be augmented at the boundary. This coupling at the boundary is best explained with the help of Figure \ref{fig:coupling}. The red particles are those of $\tau$ and the blue are those of $\alpha$. The arrows and labels represent the transitions and rates associated with the boundary coupling dynamics. For instance, if there is neither a $\tau$ or $\alpha$ particle at site 1, then at rate $\alpha$ a $\tau$-particle can enter, and at rate $\alpha'-\alpha$ an $\alpha$-particle can enter. Thus, if I am only keeping track of the $\tau$-particles, the transition from empty to occupied occurs at rate $\alpha$ and if I am tracking the $\tau+\alpha$ particles (i.e. the occurrence of their, which should match the $\tau'$ dynamic), this occurs at rate $\alpha+ (\alpha'-\alpha)=\alpha'$ as desired. Similar considerations imply that these dynamics project onto the $\tau$ and $\tau'$ dynamics marginally, and hence shows that this constitutes the desired coupling. This type of coupling is also used in the open ASEP mixing time work of Gantert, Nestoridi and Schmid \cite{gantert2020mixing}. In fact, in the case when Liggett's condition is assumed, this monotonicity seems to first have been shown in Liggett's 1975 paper \cite{Liggett75} as Corollary 3.8, based on a calculation involving the open ASEP generator.

Given the two inputs I have discussed above along with the convergence result from Section \ref{sec:conv}, let me complete the construction of the open KPZ stationary measure. As indicated above, the stochastic sandwiching and Bernoulli cases provide a way to demonstrate the H\"older bounds for $h^{(N)}(x):=N^{-1/2} h_{N}(Nx)$ (where $h_N$ is the open ASEP stationary measure in question under the Liggett condition, weak asymmetry and triple-point scaling). The calculation for the Bernoulli case is fairly simple, and is transferred to the general $u+v\neq 0$ case through the sandwiching and some simply inequalities.

Once the H\"older bounds are in place, the argument proceeds by observing that these bounds imply tightness of the initial data, and hence (through the results discussed in Section \ref{sec:conv}) also tightness of the process $h^{(N)}(\cdot,\cdot)$. Any subsequential limit will solve the open KPZ equation and any subsequential limit of the initial data $h^{(N)}(\cdot)$ will be a stationary measure (since the same held true before taking the limit). Of course, if someone proves that there is only one open KPZ stationary measure (for each pair of $u,v\in \mathbb{R}$) then this would imply that $h^{(N)}(\cdot)$ converges to it.

The above argument shows everything claimed in Section \ref{sec:const} except for the Laplace transform formula that characterizes the $u+v>0$ stationary measure. The next section explains the origin and derivation of such formulas.

\section{Matrix product ansatz and Askey-Wilson processes}\label{sec:MPA}
It is at this point in the talk version of this note that my time usually is close to cutting off (of course, much of what has been said above is severely curtailed in the hour long talk as well). So, I would like to use this space to give a sketch of the progression of ideas that form the starting point to derive the Laplace transform formula for the open KPZ stationary measure.

The first major breakthrough in studying the stationary measure for open ASEP was in Liggett's 1975 work \cite{Liggett75} where he discovered that the stationary measure satisfies a recursion relation with respect to the system size $N$. Namely, the size $N$ stationary measure is expressible in terms of the size $N-1$ case (see \cite[Theorem 3.2]{Liggett75}). In that same work, Liggett proved (based on a fairly simple generator calculation) that assuming the Liggett condition \eqref{eq:Liggett}, if $\rho_\ell=\rho_r$ then the stationary measure is Bernoulli with that density parameter. He also proved certain useful monotonicity results about the stationary measure, including with respect to changes in $N$ and in the boundary parameters. Liggett's main motivation in that work was to use open ASEP to approximate half-line and then full-line ASEP and study the convergence of different choices of initial data in that final setting to Bernoulli product measures.

A decade and a half later, there was significant renewed interest in the open ASEP stationary measure, this time from the perspective of statistical physics. The phase diagram and computation of thermodynamic quantities became centrally important. Though initially Liggett's recursion was used to make such calculations (see, e.g. \cite{DDM}), an alternative algebraic approach soon presented itself as quite useful in extracting asymptotics.

\subsection{Deriving the matrix product ansatz}
Derrida, Evans, Hakim and Pasquier introduced the {\it matrix product ansatz} (MPA) in their seminar 1993 work \cite{Derrida_1993}. They primarily focused on the case of open TASEP ($q=0$) for which they could find useful representations of the matrices. However, they also explained how to formulate the ansatz for the general open ASEP case. The MPA has found many uses since then in the study of boundary driven integrable spin chains and particle systems, and I will not try to survey this literature. Instead, let me briefly introduce the MPA for open ASEP and describe how it is derived. After that I will explain how to go from there to our open KPZ result.

As earlier in the text, let $\tau=(\tau_1,\ldots,\tau_N)\in \{0,1\}^{N}$ denote the state-space for open ASEP, and $\tau(t)$ denote the occupation process at time $t$. Then, for any $\eta,\tau \in \{0,1\}^{N}$, and $P_{\eta}(t,\tau):= \mathbb{P}(\tau(t)=\tau | \tau(0)=\eta)$, the Kolmogorov backward equation (or master equation) says  that $\partial_t P_{\eta}(t,\tau) = L^{*} P_{\eta}(t,\tau)$ where $L^*$ acts on the $\tau$ variable and represents the backwards generator of open ASEP. In particular, the stationary measure must have a zero time derivative and hence, recalling the notation $\pi_N(\tau)$ introduced earlier for the stationary measure, it must satisfy the defining relation $L^*\pi_N(\tau)=0$. Since there is no particle conservation and all states communicate in open ASEP, the Perron-Frobenius theorem implies that (up to constant scaling) there is a unique eigenfunction for $L^*$ with zero eigenvalue. This means that for any function $f_N(\tau)$ that is not identically zero and which satisfies $L^* f_N(\tau) =0$, $\pi_N(\tau)$ must equal $f_N(\tau)/Z_N$ where $Z_N$ is the sum of $f_N(\tau)$ over all $\tau \in \{0,1\}^{N}$.

Let me now write down explicitly the relation $L^* f_N(\tau)=0$, separating things out in terms of particle movement between the left reservoir and site $1$,  sites $i$ and $i+1$ in the bulk for $i=1,\ldots N-1$, and site $N$ and the right reservoir:
\begin{equation}\label{eq:genzero}
L^* f_N(\tau) :=L^*_\ell f_N(\tau) + \sum_{i=1}^{N-1} L^*_{i,i+1} f_N(\tau) + L^*_r f_N(\tau) =0
\end{equation}
where
$$
L^*_\ell f_N(\tau) := \sum_{\sigma_1} (h_\ell)_{\tau_1;\sigma_1} f_N(\sigma_1,\tau_2,\ldots \tau_N),\qquad
L^*_r f_N(\tau) := \sum_{\sigma_N} (h_r)_{\tau_N;\sigma_N} f_N(\tau_1,\ldots, \tau_{N-1},\sigma_N),
$$
and for each $i\in \{1,\ldots, N-1\}$,
$$
 L^*_{i,i+1} f_N(\tau) := \sum_{\sigma_i,\sigma_{i+1}} (h)_{\tau_i,\tau_{i+1};\sigma_i,\sigma_{i+1}} f_N(\tau_1,\ldots, \sigma_{i},\sigma_{i+1},\ldots, \tau_N).$$
In the above expressions, the sums over the $\sigma$'s are take over values in $\{0,1\}$,  $h_{\ell}$ and $h_{r}$ are $2\times 2$ matrices and $h$ is a $4\times 4$ matrix. These matrices are given as
\begin{equation}
h_{\ell}=
\left(
  \begin{array}{cc}
    -\alpha & \gamma \\
    \alpha & -\gamma \\
  \end{array}
\right),
\qquad
h=
\left(
  \begin{array}{cccc}
    0 & 0&0&0 \\
    0& -q & 1 & 0\\
    0& q & -1 & 0\\
    0 & 0&0&0 \\
  \end{array}
\right),
\qquad
h_{r}=
\left(
  \begin{array}{cc}
    -\delta & \beta \\
    \delta & -\beta \\
  \end{array}
\right)
\end{equation}
where for $h_\ell$ (and similarly $h_r$) the term $(h_{\ell})_{i;j}$ corresponds to row $i+1$ and column $j+1$ with $i,j,\in \{0,1\}$ and for $h$, the term $(h)_{i,i';j,j'}$ corresponds to row $i'+ 2i+1$ and column $j'+2j+1$ with $i,i',j,j'\in \{0,1\}$.

In order for \eqref{eq:genzero} to hold, it would be sufficient (though certainly not necessary) that there exist two constants $x_0$ and $x_1$ such that the following relations hold for all choices of $\tau_1,\ldots, \tau_N$:
\begin{align}
L^*_{\ell} f_N(\tau) &= x_{\tau_1} f_{N-1}(\hat{\tau}_1),\label{eq:telescoping1}\\
L^*_{i,i+1} f_N(\tau) &= -x_{\tau_i} f_{N-1}(\hat{\tau}_i)+x_{\tau_{i+1}} f_{N-1}(\hat{\tau}_{i+1}),\label{eq:telescoping2}\\
L^*_{r} f_N(\tau) &= -x_{\tau_N} f_{N-1}(\hat{\tau}_N),\label{eq:telescoping3}
\end{align}
where $\hat{\tau}_i$ denotes the vector $\tau$ with the $i^{th}$ coordinate (i.e., $\tau_i$) removed (e.g. $\hat{\tau}_1 = (\tau_2,\ldots, \tau_N)$).
Clearly, if these relations hold, then summing the left-hand sides yields the left-hand side of \eqref{eq:genzero} while summing the right-hand side yields 0 by telescoping.

Of course, the question is whether $f_N$ actually satisfies this relation and if so, whether there exists a manageable representation for the solution to this recursion. Inspired by earlier work in integrable systems, Derrida, Evans, Hakim and Pasquier \cite{Derrida_1993} proposed an ansatz for the form that $f_N$ could take. Consider the class of functions $\tilde{f}_N(\tau)$ that can be written in the form
$$
\tilde{f}_N(\tau_1,\ldots, \tau_N) = \langle W| (\tau_1 D + (1-\tau_1)E)\cdots (\tau_N D + (1-\tau_N)E) |V\rangle
$$
where $D$ and $E$ are matrices (possibly infinite dimensional), $\langle W|$ is a row vector and $|V\rangle$ is a column vector. Here the dimension for all of these matrices and vectors are assumed to match and the multiplication is well-defined (i.e. everything is convergent if things are infinite dimensional). For instance, if everything is one-dimensional, then the class of measures that can be defined by $\tilde{f}_N$ above (after normalizing) is exactly that of product Bernoulli measure, all with the same parameter.

Assuming the general form of $\tilde{f}_N$ above, it is possible to deduce conditions on the matrices and vectors that are necessary in order that $\tilde{f}_N$ satisfies \eqref{eq:telescoping1}-\eqref{eq:telescoping3}. Consider \eqref{eq:telescoping1}. In terms of the matrix product, this asks that
\begin{align}
 \langle W| (-\alpha E + \gamma D) \prod_{i=2}^{N} (\tau_i D + (1-\tau_i)E) |V\rangle &= x_0 \langle W|\prod_{i=2}^{N} (\tau_i D + (1-\tau_i)E) |V\rangle,\\
 \langle W| (\alpha E - \gamma D) \prod_{i=2}^{N} (\tau_i D + (1-\tau_i)E) |V\rangle &= x_1 \langle W|\prod_{i=2}^{N} (\tau_i D + (1-\tau_i)E) |V\rangle
\end{align}
where the product should be understood as ordered from left-to-right in terms of increasing index. Clearly this implies that $x_0=-x_1$ and since everything scales homogeneously, it is fine to take $x_1=1$. In order for the above relations to hold, it is sufficient then, that
$$
 \langle W| (\alpha E - \gamma D)  =  \langle W|.
$$
Similar reasoning show that in order for $\tilde{f}_N$ to satisfy \eqref{eq:telescoping3} it is sufficient that
$$
(\beta D - \delta E)|V\rangle = |V\rangle.
$$
Likewise, the relation in \eqref{eq:telescoping2} will be satisfied by $\tilde{f}_N$ as long as
$$
DE -qED = D+E.
$$
Notice that the quadratic term on the right-hand side above arises since the bulk relation \eqref{eq:telescoping2} involves summing over two particles.

To summarize, provided  non-trivial matrices and vectors (potentially infinite dimensional, though necessarily such that products are well-defined) that satisfy the three relations
\begin{equation}\label{eq.MPArelations}
DE -qED = D+E,\qquad  \langle W| (\alpha E - \gamma D)  =  \langle W|,\qquad
(\beta D - \delta E)|V\rangle = |V\rangle
\end{equation}
then necessarily the corresponding $\tilde{f}_N(\tau)$ will satisfy the defining conditions asked for $f_N(\tau)$, and hence yield the stationary measure (after normalizing). This quadratic algebra is sometimes called the DEHP algebra. Of course, there is a priori no reason to expect that there exist matrix representations to the DEHP algebra.

\subsection{Representations for the DEHP algebra}\label{sec:DEHPrep}

As an immediate application, it is now possible to use the MPA to give a sufficient condition under which the stationary measure for open ASEP is homogeneous product Bernoulli. If the matrices and vectors are all one-dimensional (which implies the Bernoulli product measure) then the DEPH algebra relations reduce to
$$de(1 -q) = d+e,\qquad  \alpha e - \gamma d  = 1 ,\qquad
\beta d - \delta e = 1
$$
for scalars $d,e$. In order for there to exist such $d$ and $e$ is suffices that the open ASEP parameters satisfies
$$
(1-q)(\alpha+\delta)(\beta+\gamma) = (\alpha+\beta+\gamma+\delta)(\alpha\beta - \gamma \delta).
$$
If we further assume that Liggett's condition holds (thus expressing $\gamma$ in terms of $\alpha$, and $\delta$ in terms of $\beta$) then the above equation is satisfied when $\alpha +\beta=1$, which was precisely the condition mentioned in Section \ref{sec:asym}.

Putting aside the one-dimensional case, \cite{Derrida_1993} showed that any other matrix representation to \eqref{eq.MPArelations} must be infinite dimensional. Finding such representations is a highly non-trivial task.  When $q=0$, \cite{Derrida_1993} provided a few such representations. Different matrix representations have proven useful in making various types of calculations involving the large $N$ limit of the stationary measure, for example including computing the stationary current, and other correlation functions.

For the general $q\neq 0$ and general $\alpha,\beta,\gamma,\delta$ parameter case, it took about a decade until Uchiyama, Sasamoto and Wadati \cite{USW} provided the first representations to the DEHP algebra. Their infinite dimensional matrix representations (that I will call the USW representation) were written in terms of the Askey-Wilson orthogonal polynomials (in particular, the associated Jacobi matrix). This remarkable link between orthogonal polynomials and the open ASEP stationary measure was already present in a simpler case (when $\gamma=\delta=0$) in earlier work of just Sasamoto \cite{Sasamoto_1999} in which the Al-Salam-Chihara polynomials replaced the Askey-Wilson polynomials. Another notable development in this area came in work of Corteel and Williams \cite{corteel2011} who, using a different matrix representation, found a combinatorial description for the open ASEP stationary measure in terms of certain types of tableaux combinatorics. A nice exposition on this direction and further developments around it can be found in Williams' recent expository piece \cite{WilliamsLMS}.

Since my aim here is explain how the work of \cite{USW} leads to the Askey-Wilson process formulas of Bryc and Weso{\l}owski \cite{BrycWes17}, I will recall a variant of the USW matrix representation for the DEHP algebra in the notation of  \cite{BrycWes17}.

It will be helpful to work with the parameterizations of $\alpha,\beta,\gamma,\delta$ in terms of $A,B,C,D$ as facilitated by the bijection in
\eqref{eq_ABCD}. The  USW matrix representation is infinite dimensional and written in terms of one-sided infinite matrices and vectors (i.e., indexed by natural numbers). I will use bold face variables to distinguish this representation as well as associated matrices used in describing it. The row vector $\langle \mathbf{W}|=[1,0,0,\ldots]$ and column vector $|\mathbf{V}\rangle = [1,0,0,\ldots]^T$ are both simple. The matrices $\mathbf{E}$ and $\mathbf{D}$ are more complicated and can be written in terms of the identity matrix $\mathbf{I}$ and two other matrices $\mathbf{x}$ and $\mathbf{y}$ as
\begin{equation}\label{eqDE}
\mathbf{D} = \frac{1}{1-q} \mathbf{I} + \frac{1}{\sqrt{1-q}}\mathbf{x},\qquad \mathbf{E} = \frac{1}{1-q} \mathbf{I} + \frac{1}{\sqrt{1-q}}\mathbf{y}.
\end{equation}
The  $\mathbf{x}$ and $\mathbf{y}$ matrices are tridiagonal and admit explicit formulas for the coefficients. Since I do not want to get lost in the (very important) details here, I will abstain from recording them precisely, but rather just focus on their key properties. First and foremost, they are such that $\mathbf{D}$ and $\mathbf{E}$, along with the simple $\langle \mathbf{W}|$ and $|\mathbf{V}\rangle$, satisfy the DEHP algebra. That they satisfy the quadratic relation in the DEHP algebra is equivalent to the $q$-commutation relation $\mathbf{x}\mathbf{y} - q\mathbf{y}\mathbf{x} = \mathbf{I}$. This relation, in conjunction with the knowledge of the first entry of the three non-trivial diagonals in both $\mathbf{x}$ and $\mathbf{y}$ uniquely determine their values. The values of these non-trivial entries can be determined from the two boundary relations in the DEHP algebra.

\subsection{Askey-Wilson polynomials and processes}\label{sec:AWPP}

The other key property of the $\mathbf{x}$ and $\mathbf{y}$ matrices is that they appear in the Jacobi matric that describes the three step recurrence for Askey-Wilson orthogonal polynomials. Let me briefly introduce these polynomials since they will be needed in order to relate the MPA to the Askey-Wilson process (which is not defined in terms of these polynomials, but rather the orthogonality measure associated to them).

I will use the notation of \cite{BrycWes17} whereby they write the Askey-Wilson polynomials as $\bar{w}_n(x;a,b,c,d,q)$. The parameters should be assumed to satisfy  $q\in (-1,1)$, $a,b,c,d\in \mathbb{C}$ with $abcd,abcdq,ab,abq\notin [1,\infty)$. These polynomials are typically described in terms of their three term recurrence, though I will not write this down here explicitly. Under some additional conditions on parameters (that I will also suppress here) these polynomials will be orthogonal with respect to a probability distribution known as the {\it Askey-Wilson measure}. This measure has the form
\begin{equation}\label{eq:AWmeasure}
\nu(dx;a,b,c,d,q) := f(x;a,b,c,d,q) \mathbf{1}_{|x|<1} + \sum_{y\in F(a,b,c,d,q)} p(y;a,b,c,d,q) \delta_{y}(dx).
\end{equation}
The first part above is an absolutely continuous measure on $x\in(-1,1)$ with density
\begin{equation}\label{eq:AWden}
f(x;a,b,c,d,q) = \frac{(q,ab,ac,ad,bc,bd,cd;q)_{\infty}}{2\pi (abcd;q)_{\infty} \sqrt{1-x^2}} \left|\frac{(e^{2\iota \theta};q)_{\infty}}{(ae^{\iota \theta},be^{\iota \theta},ce^{\iota \theta},de^{\iota \theta};q)_{\infty}}\right|^2
\end{equation}
where we recall $(a;q)_{\infty} = (1-a)(1-qa)(1-q^2a)\cdots$, $(a_1,\ldots, a_k;q)_{\infty} = (a_1;q)_{\infty}\cdots (a_k;q)_{\infty}$, and $\theta$ is defined by the relation $\cos(\theta) =x$. The second part of \eqref{eq:AWmeasure} is an atomic measure supported at locations $y$ in the finite set $F(a,b,c,d,q)$ with masses $p(y;a,b,c,d,q)$. The points and their masses are explicit as well, though I will not record them here. Though I will not explain it, these atoms only become important when the open KPZ parameters $u,v$ are not both positive. The formula \eqref{eq:laplace} that I gave earlier for the Laplace transform of the stationary open KPZ height increment from 0 to 1 was exactly in the case of $u,v>0$ in which case it is only the absolutely continuous part above in \eqref{eq:AWmeasure} that is present.

In a spirit similar (albeit considerably more involved)  to how one uses the Gaussian distribution (which is the orthogonality measure for Hermite polynomials) to define Brownian motion, Bryc and Weso{\l}owski \cite{BrycWes10} used the measure $\nu$ to define a process that they called the {\it Askey-Wilson process}. The existence and path properties (e.g. continuity versus jumps) of this process are non-trivial. It is defined only for an interval $I$ of time that is dependent upon the five parameters $A,B,C,D,q$ that determine it. These parameters are also subject to certain conditions for existence of the process. Following \cite{BrycWes10} and \cite{BrycWes17}, let $Y = (Y_t)_{t\in I}$ denote this process. It is a Markov process with marginal distribution for each $t\in I$ given by the measure
$$
\pi_t(dx) := \nu(dx;A\sqrt{t},B\sqrt{t},C/\sqrt{t},D/\sqrt{t},q)
$$
and with transition probabilities for $s<t$ both in $I$ of the form
$$
P_{s,t}(x,dy) := \nu(dy; A\sqrt{t},B\sqrt{t}, \sqrt{s/t} (x+\sqrt{x^2-1}),\sqrt{s/t}(x-\sqrt{x^2-1}),q).
$$
It takes some work to show that this is well-defined, i.e., satisfying Chapman-Kolmogorov.

Bryc and Weso{\l}owski's original motivation for introducing the Askey-Wilson process came from the study of {\it quadratic harness} which are stochastic  processes $X_t$ such that the expectation of $X_t$ and of $X_t^2$, conditioned on the process outside an interval containing $t$, are given as linear and quadratic functions (respectively) of the values of the process $X$  at the boundary of said interval. Brownian motion or a Poisson jump process are arguably the simplest examples of such harnesses.  The question addressed in \cite{BrycWes10} was how to characterize the space of all standard quadratic harnesses (standard means that $X_t$ has the same mean and covariance as standard Brownian motion), and it turned out that this was achieved by the Askey-Wilson processes.

The first time I heard about Askey-Wilson processes was at a conference on ``Orthogonal polynomials, applications in statistics and stochastic processes'' in 2010 at Warwick. This was quite a notable event for me since it was among the first times I gave a talk at a conference, especially one overseas. I have a vivid memory of sitting in the lecture theater and listening to Weso{\l}owski give a talk entitled ``Quadratic harnesses and Askey-Wilson polynomials''. One of the reasons it was so vivid is that I was completely confused as to why anyone would study these seemingly very complicated processes. I suppose the moral of this brief anecdote is that when someone gives a talk about something that seems very  complicated, it is often worth paying attention even if you do not immediately understand why.

A property of the Askey-Wilson polynomials and process is that when combined, they form a family of orthogonal martingales. Specifically, let
$$
Z_t := \frac{2\sqrt{t}}{\sqrt{1-q}} Y_t
$$
and define the (infinite) row vector-valued function $(x,t)\mapsto \langle \mathbf{r}_t(x)|$ where
$$
\langle \mathbf{r}_t(x)| := [ r_0(x;t),r_1(x;t),\ldots]
$$
and where $r_n(x;t)$ are polynomials of degree $n$ in the variable $x$ and given in terms of the Askey-Wilson polynomial $\bar{w}_n$ by
$$
r_n(x;t) = t^{n/2} \bar{w}_n\left(\frac{\sqrt{1-q}}{2\sqrt{2}}; A\sqrt{t},B\sqrt{t},C/\sqrt{t},D/\sqrt{t},q\right).
$$
The process $\langle \mathbf{r}_t(Z_t)|$ satisfies the following properties (proved in \cite{BrycWes10}):
\begin{enumerate}
\item $r_0(x;t)=1$ for all $x$,
\item $\mathbb{E}[r_n(Z_t;t)r_m(Z_t;t)] = 0$ if $m\neq n$,
\item $\mathbb{E}[r_n(Z_t;t) | Z_s] = \mathbb{E}[r_n(Z_t;t) | \mathcal{F}_s] = r_n(Z_s;s)$ for $s\leq t$ and $\mathcal{F}_s=\sigma(Z_v:v\leq s)$ the sigma-algebra generated by $Z_v$ for all $v\leq s$.
\end{enumerate}

The three term recurrence relation for Askey-Wilson polynomials can be rewritten in terms of the $r_n$ version of these polynomials in the following succinct manner:
\begin{equation}\label{eq:AWthreeterm}
x\langle \mathbf{r}_t(x)| = \langle \mathbf{r}_t(x)| (t \mathbf{x} + \mathbf{y})
\end{equation}
where $\mathbf{x}$ and $\mathbf{y}$ are the matrices discussed earlier around \eqref{eqDE}. Since these are tridiagonal matrices, they imply that $xr_n(x;t)$ can be expressed as a sum of $r_{n-1}(x;t)$, $r_{n}(x;t)$ and $r_{n+1}(x;t)$ with coefficients that depend on $t$ in a linear manner but not on the variable $x$.

\subsection{Connecting matrix product ansatz and Askey-Wilson processes}

With all of the notation and properties introduced above, I am now in a position to relate the beautiful proof that appears in \cite{BrycWes17} and relates the open ASEP stationary measure, via the matrix product ansatz and the Uchiyama-Sasamoto-Wadati representation, to the Askey-Wilson process.

The matrix product ansatz gives a formula for the stationary measure. However, by summing over all states it also give a compact representation for the following generating function for $\pi_N(\tau)$ (recall the notation on the left-hand side means averaging a function of $\tau$ against this measure)
$$
\left\langle \prod_{j=1}^{N} t_j^{\tau_j}\right\rangle_N = \frac{\langle \mathbf{W}| \prod_{j=1}^{N} (\mathbf{E} + t_j \mathbf{D})|\mathbf{V}\rangle}{\langle \mathbf{W}| \prod_{j=1}^{N} (\mathbf{E} + \mathbf{D})|\mathbf{V}\rangle}.
$$
As before, the product is ordered from left to right in increasing order of indices. Focusing on the numerator and using the USW representation \eqref{eqDE} of the DEHP algebra yields
$$
\langle \mathbf{W}| \prod_{j=1}^{N} (\mathbf{E} + t_j \mathbf{D})|\mathbf{V}\rangle = (1-q)^{-N} \langle \mathbf{W}| \prod_{j=1}^{N} \big((1+t_j)\mathbf{I} + \sqrt{1-q}(t_j\mathbf{x}+\mathbf{y})\big)|\mathbf{V}\rangle : = (1-q)^{-N} \Pi.
$$
Since $r_0(x;t)\equiv 1$ and since $\mathbb{E}[r_n(Z_t;t)]=\mathbb{E}[r_n(Z_t;t)r_0(Z_t;t)]=0$ if $n\neq 0$, the vector $\langle \mathbf{W}|$ can be rewritten in terms of the vector $\mathbb{E}[\langle \mathbf{r}_t(Z_t)|]$ for any choice of $t$. By linearity of the expectation,
$$
\Pi = \mathbb{E}\left[ \langle \mathbf{r}_{t_1}(Z_{t_1})| \prod_{j=1}^{N}  \big((1+t_j)\mathbf{I} + \sqrt{1-q}(t_j\mathbf{x}+\mathbf{y})\big) |\mathbf{V}\rangle \right].
$$
The three term recurrence relation \eqref{eq:AWthreeterm} can now be applied to peal off the first term in the product above. In particular, observe that for any $t$,
$$
 \langle \mathbf{r}_{t}(Z_{t})| \big((1+t)\mathbf{I} + \sqrt{1-q}(t\mathbf{x}+\mathbf{y})\big) =  (1+t + \sqrt{1-q}Z_t)\langle \mathbf{r}_t(Z_t)|.
$$
The term involving $\mathbf{I}$ follows since that is just the identity matrix, and the other follows from \eqref{eq:AWthreeterm}. This means that
$\Pi$ can be rewritten as
$$
\Pi = \mathbb{E}\left[  (1+t + \sqrt{1-q}Z_{t_1})  \langle \mathbf{r}_{t_1}(Z_{t_1})| \prod_{j=2}^{N}  \big((1+t_j)\mathbf{I} + \sqrt{1-q}(t_j\mathbf{x}+\mathbf{y})\big) |\mathbf{V}\rangle \right].
$$
Provided that $t_2\geq t_1$, the martingale property for each $r_n(Z_t;t)$ implies that
$$
 \langle \mathbf{r}_{t_1}(Z_{t_1})| \prod_{j=2}^{N}  \big((1+t_j)\mathbf{I} + \sqrt{1-q}(t_j\mathbf{x}+\mathbf{y})\big) |\mathbf{V}\rangle
 =
 \mathbb{E}\left[\langle \mathbf{r}_{t_2}(Z_{t_2})| \prod_{j=2}^{N} \big((1+t_j)\mathbf{I} + \sqrt{1-q}(t_j\mathbf{x}+\mathbf{y})\big) |\mathbf{V}\rangle | \mathcal{F}_{t_1}\right].
$$
Plugging this into the above expression for $\Pi$ and using the tower property for martingales to remove the conditional expectation yields
$$
\Pi  = \mathbb{E}\left[  (1+t + \sqrt{1-q}Z_{t_1})
 \langle \mathbf{r}_{t_2}(Z_{t_2})| \prod_{j=2}^{N}  \big((1+t_j)\mathbf{I} + \sqrt{1-q}(t_j\mathbf{x}+\mathbf{y})\big) |\mathbf{V}\rangle | \right].
$$
This procedure can be repeated providing $t_1\leq t_2\leq t_3\leq \cdots \leq t_N$, yielding a final result of
$$
\Pi =  \mathbb{E}\left[  \prod_{j=1}^{N} (1+t_j + \sqrt{1-q}Z_{t_j}) \langle \mathbf{r}_{t_N}(Z_{t_N})| \mathbf{V}\rangle \right].
$$
However, since $|\mathbf{V}\rangle = [1,0,\ldots]^T$ and $\langle \mathbf{r}_{t_N}(Z_{t_N})| \mathbf{V}\rangle = r_0(Z_{t_N};t_N) = 1$. Thus,
$$
\Pi =  \mathbb{E}\left[  \prod_{j=1}^{N} (1+t_j + \sqrt{1-q}Z_{t_j}) \right].
$$
In the above calculation it was important that the Askey-Wilson process be defined for the values of $t_j$ used, and for the parameters $A,B,C$ and $D$ that encode the boundary parameters for open ASEP. Without going into the details, provided that $AC<1$, the Askey-Wilson process is well-defined for $t\in (0,\infty)$ and thus the above argument proves the following result  \cite[Theorem 1]{BrycWes10}: For $0<t_1\leq \cdots \leq t_N$,
\begin{equation}\label{eq:ASEPAWP}
\left\langle \prod_{j=1}^{N} t_j^{\tau_j}\right\rangle_N  = \frac{ \mathbb{E}\left[  \prod_{j=1}^{N} (1+t_j + \sqrt{1-q}Z_{t_j}) \right]}{ \mathbb{E}\left[  (2 + \sqrt{1-q}Z_{1})^N \right]}
\end{equation}
where the $\langle \cdot\rangle_N$ on the left-hand side is the stationary state expectation for open ASEP and the $\mathbb{E}[\cdot]$ on the right-hand side is the expectation with respect to the Askey-Wilson process (with the two processes related in terms of their parameters $A,B,C,D$ and $q$).  The condition $AC<1$ corresponds exactly to the fan region discussed in Section \ref{sec:phasediagram} and ultimately is the reason why the open KPZ stationary measure Laplace transform is only computed for $u+v>0$ (the limit of the fan region).

\subsection{Coming full circle to asymptotics}\label{sec:fullcircle}

To close, I would like to briefly explain why \eqref{eq:ASEPAWP} is quite useful for asymptotics. Recall from earlier the formula \eqref{eq:netchange} for the Laplace transform for the total change in height for the open KPZ equation across the interval $[0,1]$. The microscopic analog for this is $h_N(N)=2(\tau_1+\cdots +\tau_N)-N$.

If $t_1=\cdots =t_N=e^{2s}$ for $s\in \mathbb{R}$ then \eqref{eq:ASEPAWP} implies the Laplace transform formula
\begin{equation}\label{eq:shn}
\langle e^{s h_N(N)}\rangle  = e^{-sN} \frac{ \mathbb{E}\left[   (1+e^{2s}+ \sqrt{1-q}Z_{e^{2s}})^N \right]}{ \mathbb{E}\left[ (2 + \sqrt{1-q}Z_{1})^N\right]}.
\end{equation}

The right-hand side numerator (and similarly denominator) can be written as an integral against the law of $Z_{e^{2s}}$, i.e., against the (scaled) Askey-Wilson orthogonality measure. As $N$ grows, the complexity of the right-hand side does not. In fact, as is often the case in asymptotic analysis, owing to the power of $N$, the right-hand side will actually simplify in the $N\to\infty$ limit, even when the special scalings from Section \ref{sec:aseptokpz} are applied. This should be compared to the original form of the matrix product ansatz which requires multiplication of $N$ matrices, the complexity of which grows considerably, even if they are tridiagonal.

I should be clear that the asymptotics of a formula like \eqref{eq:shn} is not simple. However, it can be done. In \cite{BrycWes17}, they apply their Askey-Wilson process formula \eqref{eq:ASEPAWP} to prove a large deviation principle for the stationary measure of open ASEP with fixed parameters as $N$ goes to infinity. Soon after, Bryc and Wang \cite{BW} used \eqref{eq:ASEPAWP} to study the fluctuation scaling limit for the stationary measure, again for all parameters fixed and $N\to \infty$.


My results with Knizel, in particular, the open KPZ Laplace transform formula discussed in Section \ref{sec:const} also proceed through \eqref{eq:ASEPAWP}. In our case, all of the parameters are being scaled in an $N^{-1/2}$ window around their limiting values:  $q$ approaching $1$, $\alpha,\beta,\gamma,\delta$ approaching $1/2$ and (in order to accommodate the scaling of $h_N(N)$),  $s$ approaching $0$. Inputting these scalings into a formula like \eqref{eq:shn}, it eventually becomes clear that the main contribution to the integral against the Askey-Wilson process marginal distribution comes when $Z$ is within order $N^{-1}$ of $1$. Zooming into this scale eventually (in the $N\to\infty$ limit) leads to formulas involving a {\it tangent process} to the Askey-Wilson process -- the continuous dual Hahn process mentioned earlier in  Section \ref{sec:const}. In terms of \eqref{eq:shn}, all of this scaling also introduces a diverging Jacobian factor which is compensated by the decay of the integrand. The point-wise limit of the integrand and measure exists based on convergence of $q$-gamma functions to gamma functions (recall that $q$ is tending to $1$ and there are lots of $q$-Pochhammer symbols). However, in order to conclude that the integral itself converges (and hence deduce formulas like \eqref{eq:netchange}) requires an application of the dominated convergence theorem. This, in turn, relies on uniform control over the behavior of $q$-gamma functions with its variable varying in vertical strips in the complex plane of height of order $N$ (though $Z$ has bounded support, after scaling by order $N^{-1}$ around $Z=1$, the support becomes of order $N$). So, to bring things full circle, it is exactly in proving this type of dominated convergence bounds that the results of Section \eqref{Sec:qaside} become necessary.

\bibliographystyle{alpha}
\bibliography{refs}

\begin{thebibliography}{BKWW21}

\bibitem[AW85]{AskeyWilson}
R.~Askey and J.~Wilson.
\newblock {Some basic hypergeometric orthogonal polynomials that generalize
  Jacobi polynomials}.
\newblock {\em Memoirs of the AMS}, 54, 1985.

\bibitem[BC95]{BertiniCancrini}
L.~Bertini and N.~Cancrini.
\newblock {The stochastic heat equation: Feynman-Kac formula and
  intermittence}.
\newblock {\em Journal of Statistical Physics}, 78:1377--1401, 1995.

\bibitem[BCK14]{BCKFullLine}
Y.~Bakhtin, E.~Cator, and K.~Khanin.
\newblock {Space-time stationary solution for the Burgers equation}.
\newblock {\em Journal of the AMS}, 27(1):193--238, 2014.

\bibitem[BE07]{Blythe_2007}
R.~A. Blythe and M.~R. Evans.
\newblock Nonequilibrium steady states of matrix-product form: a solver's
  guide.
\newblock {\em Journal of Physics A: Mathematical and Theoretical},
  40(46):R333--R441, oct 2007.

\bibitem[BG97]{cmp/1158328658}
L.~Bertini and G.~Giacomin.
\newblock {Stochastic Burgers and KPZ equations from particle systems}.
\newblock {\em Communications in Mathematical Physics}, 183(3):571 -- 607,
  1997.

\bibitem[BK21]{bryc2021markov}
W.~Bryc and A.~Kuznetsov.
\newblock {Markov limits of steady states of the KPZ equation on an interval},
  2021.
\newblock arXiv:2109.04462.

\bibitem[BKWW21]{BKWW}
W.~Bryc, A.~Kuznetsov, Y.~Wang, and J.~Weso{\l}owski.
\newblock {Markov processes related to the stationary measure for the open KPZ
  equation}, 2021.
\newblock arXiv:2105.03946.

\bibitem[BLD21]{BLD}
G.~Barraquand and P.~Le~Doussal.
\newblock {Steady state of the KPZ equation on an interval and Liouville
  quantum mechanics}, 2021.
\newblock arXiv:2105.15178.

\bibitem[Bry22]{BRYC2022185}
W.~Bryc.
\newblock On the continuous dual hahn process.
\newblock {\em Stochastic Processes and their Applications}, 143:185--206,
  2022.

\bibitem[BW10]{BrycWes10}
W.~Bryc and J.~Weso{\l}owski.
\newblock {A}skey--{W}ilson polynomials, quadratic harnesses and martingales.
\newblock {\em Annals of Probability}, 38:1221--1262, 2010.

\bibitem[BW17]{BrycWes17}
W.~Bryc and J.~Weso{\l}owski.
\newblock {A}symmetric {S}imple {E}xclusion {P}rocess with open boundaries and
  quadratic harnesses.
\newblock {\em Journal of Statistical Physics}, 167:383--415, 2017.

\bibitem[BW19]{BW}
W.~Bryc and Y.~Wang.
\newblock {Limit fluctuations for density of {A}symmetric {S}imple {E}xclusion
  {P}rocesses with open boundaries}.
\newblock {\em Annals of the Institute Henri Poincar\'e B}, 55:2169--2194,
  2019.

\bibitem[CH16]{CHammond}
I.~Corwin and A.~Hammond.
\newblock {KPZ line ensemble}.
\newblock {\em Probability Thoery and Related Fields}, 166:67--185, 2016.

\bibitem[CK21]{CK}
Ivan Corwin and Alisa Knizel.
\newblock {Stationary measure for the open KPZ equation}.
\newblock 2021.
\newblock arXiv:2103.12253.

\bibitem[CMY98]{comtet_monthus_yor_1998}
A.~Comtet, C.~Monthus, and M.~Yor.
\newblock {Exponential functionals of Brownian motion and disordered systems}.
\newblock {\em Journal of Applied Probability}, 35(2):255--271, 1998.

\bibitem[Cor12]{doi:10.1142/S2010326311300014}
I.~Corwin.
\newblock {The Kardar-Parisi-Zhang equation and university class}.
\newblock {\em Random Matrices: Theory and Applications}, 01(01):1130001, 2012.

\bibitem[Cor16]{CorwinAMS}
I.~Corwin.
\newblock {Kardar-Parisi-Zhang universality}.
\newblock {\em Notices of the AMS}, 63:230--239, 2016.

\bibitem[Cor18]{ICAMS}
I.~Corwin.
\newblock {Exactly solving the KPZ equation}.
\newblock {\em In ``Random Growth Models'', Proceedings of Symposia in Applied
  Mathematics. Editors: M. Damron, F. Rassoul-Aghna and T. Sepp\"al\"ainen},
  75, 2018.

\bibitem[CS18]{CS}
I.~Corwin and H.~Shen.
\newblock Open $\text{ASEP}$ in the weakly asymmetric regime.
\newblock {\em Communications on Pure and Applied Mathematics}, 71:2065--2128,
  2018.

\bibitem[CS20]{CorwinShenBAMS}
I.~Corwin and H.~Shen.
\newblock Some recent progress in singular stochastic partial differential
  equations.
\newblock {\em Bulletins of the AMS}, 57:409--454, 2020.

\bibitem[CW11]{corteel2011}
S.~Corteel and L.~K. Williams.
\newblock {Tableaux combinatorics for the asymmetric exclusion process and
  Askey-Wilson polynomials}.
\newblock {\em Duke Mathematical Journal}, 159(3):385--415, 09 2011.

\bibitem[Daa94]{Daalhuis}
A.~B.~O. Daalhuis.
\newblock Asymptotic expansions for q-gamma, q-exponential, and q-bessel
  functions.
\newblock {\em J. Math. Anal. Appl.}, 186:896--913, 1994.

\bibitem[DDM92]{DDM}
B.~Derrida, E.~Domany, and D~Mukamel.
\newblock An exact solution of a one-dimensional asymmetric exclusion model
  with open boundaries.
\newblock {\em Journal of Statistical Physics}, 69:667--687, 1992.

\bibitem[DEHP93]{Derrida_1993}
B.~Derrida, M.~R. Evans, V.~Hakim, and V.~Pasquier.
\newblock Exact solution of a 1d asymmetric exclusion model using a matrix
  formulation.
\newblock {\em Journal of Physics A}, 26(7):1493--1517, 1993.

\bibitem[DEL04]{DEL}
B.~Derrida, C.~Enaud, and J.~L. Lebowitz.
\newblock {The asymmetric exclusion process and Brownian excursions}.
\newblock {\em Journal of Statistical Physics}, 115:365--382, 2004.

\bibitem[DGR19]{alex2019stationary}
A.~Dunlap, C.~Graham, and L.~Ryzhik.
\newblock Stationary solutions to the stochastic burgers equation on the line,
  2019.
\newblock arXiv:1910.07464.

\bibitem[DT16]{Dembo2016}
A.~Dembo and L.-C. Tsai.
\newblock {Weakly Asymmetric Non-Simple Exclusion Process and the
  Kardar--Parisi--Zhang Equation}.
\newblock {\em Communications in Mathematical Physics}, 341:219--261, 2016.

\bibitem[EKMS]{EKMS}
W.~E, K.~Khanin, A.~Mazel, and Ya. Sinai.
\newblock Invariant measures for burgers equation with stochastic forcing.
\newblock {\em Annals of Mathematics}, 151:877--960.

\bibitem[FQ15]{FQ14}
T.~{Funaki} and J.~{Quastel}.
\newblock {KPZ equation, its renormalization and invariant measures}.
\newblock {\em Stochastic Partial Differential Equations: Analysis and
  Computation}, 3:159--220, 2015.

\bibitem[Gan21]{Ganguly}
S.~Ganguly.
\newblock {Random metric geometries on the plane and Kardar-Parisi-Zhang
  universality}.
\newblock {\em Notices of the AMS}, 69:26--35, 2021.

\bibitem[Gar88]{Gartner}
J.~Gartner.
\newblock {Convergence towards Burgers' equation and propagation of chaos for
  weakly asymmetric exclusion processes}.
\newblock {\em Stochastic Processes and Applications}, 27:233--260, 1988.

\bibitem[GH19]{GH19}
M.~Gerencs\'er and M.~Hairer.
\newblock {Singular SPDEs in domains with boundaries}.
\newblock {\em Probability Theory and Related Fields}, 173:697--758, 2019.

\bibitem[GLM17]{10.1214/16-AAP1200}
P.~Gon\c{c}alves, C.~Landim, and A.~Milan\'es.
\newblock {Nonequilibrium fluctuations of one-dimensional boundary driven
  weakly asymmetric exclusion processes}.
\newblock {\em The Annals of Applied Probability}, 27:140 -- 177, 2017.

\bibitem[GNS20]{gantert2020mixing}
N.~Gantert, E.~Nestoridi, and D.~Schmid.
\newblock Mixing times for the simple exclusion process with open boundaries,
  2020.
\newblock arXiv:2003.03781.

\bibitem[GP18]{gubinelli2018infinitesimal}
M.~Gubinelli and N.~Perkowski.
\newblock {The infinitesimal generator of the stochastic Burgers equation},
  2018.
\newblock arXiv:1810.12014.

\bibitem[HHT15]{HHTak}
T.~Halpin-Healy and K.~Takeuchi.
\newblock {A KPZ cocktail -- shaken, not stirred: Toasting 30 years of
  kinetically roughened surfaces}.
\newblock {\em Journal of Statistical Physics}, 160:794--814, 2015.

\bibitem[HM18]{hairer2018}
M.~Hairer and J.~Mattingly.
\newblock {The strong Feller property for singular stochastic PDEs}.
\newblock {\em Annals of the Institute Henri Poincar\'e B}, 54:1314--1340,
  2018.

\bibitem[HY04]{HariyaYor}
Y.~Hariya and M.~Yor.
\newblock {Limiting distributions associated with moments of exponential
  Brownian functionals}.
\newblock {\em Studia Scientiarum Mathematicarum Hungarica}, 41:193--242, 2004.

\bibitem[Kat03]{852f802e715846339a94147fc8f7f3bb}
Masanori Katsurada.
\newblock Asymptotic expansions of certain q-series and a formula of ramanujan
  for specific values of the riemann zeta function.
\newblock {\em Acta Arithmetica}, 107(3):269--298, 2003.

\bibitem[Lig75]{Liggett75}
T.~Liggett.
\newblock Ergodic theorems for the asymmetric simple exclusion process.
\newblock {\em Transactions of the American Mathematical Society},
  213:237--261, 1975.

\bibitem[Lig77]{Liggett77}
T.~Liggett.
\newblock Ergodic theorems for the asymmetric simple exclusion process ii.
\newblock {\em Annals of Probability}, 5(5):795--801, 1977.

\bibitem[Mci99]{McIntosh}
R.~J. Mcintosh.
\newblock Some asymptotic formulae for q-shifted factorials.
\newblock {\em Ramanujan Journal}, 3:205--214, 1999.

\bibitem[MGP68]{MGP}
C.~T. Macdonald, J.~H. Gibbs, and A.~C. Pipkin.
\newblock {Kinetics of biopolymerization on nucleic acid templates}.
\newblock {\em Biopolymers}, 6:1--5, 1968.

\bibitem[Moa84]{Moak}
D.~S. Moak.
\newblock {The q-analogue of Stirling's formula}.
\newblock {\em Rocky Mountain Journal of Mathematics}, 14:403--413, 1984.

\bibitem[MY05]{10.1214/154957805100000159}
H.~Matsumoto and M.~Yor.
\newblock {Exponential functionals of Brownian motion, I: Probability laws at
  fixed time}.
\newblock {\em Probability Surveys}, 2:312 -- 347, 2005.

\bibitem[Par19]{Parekh}
S.~Parekh.
\newblock {The KPZ Limit of ASEP with Boundary}.
\newblock {\em Communications in Mathematical Physics}, 365:569--649, 2019.

\bibitem[QS15]{Quastelspohn}
J.~Quastel and H.~Spohn.
\newblock {The one-dimensional KPZ equation and its universality class}.
\newblock {\em Journal of Statistical Physics}, 160:965--984, 2015.

\bibitem[Ros19]{rosati2019synchronization}
T.~C. Rosati.
\newblock {Synchronization for KPZ}, 2019.
\newblock arXiv:1907.06278.

\bibitem[Sas99]{Sasamoto_1999}
Tomohiro Sasamoto.
\newblock One-dimensional partially asymmetric simple exclusion process with
  open boundaries: orthogonal polynomials approach.
\newblock {\em Journal of Physics A: Mathematical and General},
  32(41):7109--7131, sep 1999.

\bibitem[Sch21]{Schmid}
D.~Schmid.
\newblock Mixing times for the tasep in the maximal current phase, 2021.
\newblock arXiv:2104.12745.

\bibitem[Spi70]{SPITZER1970246}
Frank Spitzer.
\newblock Interaction of markov processes.
\newblock {\em Advances in Mathematics}, 5(2):246--290, 1970.

\bibitem[USW04]{USW}
M.~Uchiyama, T.~Sasamoto, and M.~Wadati.
\newblock {Asymmetric simple exclusion process with open boundaries and
  Askey{\textendash}Wilson polynomials}.
\newblock {\em Journal of Physics A}, 37(18):4985--5002, 2004.

\bibitem[Wil22]{WilliamsLMS}
L.~K. Williams.
\newblock {The combinatorics of hopping particles and positivity in Markov
  chains}.
\newblock 2022.
\newblock arXiv:2202.00214.

\bibitem[Yan21]{yang2021kardarparisizhang}
K.~Yang.
\newblock {Kardar-Parisi-Zhang Equation from Non-Simple Variations on
  Open-ASEP}, 2021.

\bibitem[Zha14]{Z}
R.~Zhang.
\newblock On asymptotics of the $q$-exponential and $q$-gamma functions.
\newblock {\em Journal of Mathematical Analysis and Applications},
  411:522--529, 2014.

\end{thebibliography}

\end{document}